\documentclass[11pt,dvips]{amsart}  
\usepackage{amssymb,amsmath,amsfonts}  
\usepackage[dvips]{graphicx,color}
\usepackage{cite}
\allowdisplaybreaks[4]
\definecolor{pgray}{gray}{0.8}

\newtheorem{theorem}{Theorem}[section]

\newtheorem{proposition}[theorem]{Proposition}

\numberwithin{equation}{section}

\title{\mbox{}}

\begin{document}
\begin{center}
{\bf \LARGE{
	Nonlinear distortion of symmetry in solutions to the convection-diffusion equation of Burgers type
}}\\
\vspace{5mm}
Masakazu Yamamoto\footnote{e-mail : \texttt{mk-yamamoto@gunma-u.ac.jp}}\\Graduate School of Science and Technology, Gunma University
\end{center}
\maketitle
\vspace{-15mm}
%
\begin{abstract}
In this paper, the initial value problem of the convection-diffusion equation of Burgers type is treated.
In the asymptotic profile of solutions, the nonlinearity of the equation is reflected.
Regarding the solutions to this model, the Spanish school in the 1990s performed asymptotic expansions based on the linear diffusion.
Those profiles exhibit symmetries characteristic of linear phenomena.
In this paper, the distortion of symmetry arising from the nonlinear effects is described explicitly.
Furthermore, it is demonstrated that the extent of this distortion differs significantly depending on the parity of the spatial dimension.
This contradicts the conventional expectation that the manifestation of nonlinearity depends on the scale of the equation.
This interpretation is supported by comparison with similar Navier--Stokes equations.
The Burgers type is applicable as an indicator for considering several bilinear problems.
\end{abstract}

\section{Introduction}
We study the following initial-value problem for $n \ge 2$:
\begin{equation}\label{cd}
	\left\{\begin{array}{lr}
		\partial_t u - \Delta u = \boldsymbol{a} \cdot \nabla (u^2),	&	t>0,~ x \in \mathbb{R}^n,\\
		u(0,x) = u_0 (x),								&	x \in \mathbb{R}^n,
	\end{array}\right.
\end{equation}
where $\boldsymbol{a} \in \mathbb{R}^n$ is a constant and $u_0$ is a given initial data.
Especially, the cases  $n = 2,3$ and $4$ are focused.
This equation is the bilinear convection-diffusion equation for positive density, or the Burgers equation in higher dimensions.
However, the initial data does not need to be positive.
The reason for dealing with such a limited model will be explained after presenting our results.
In any cases, our equation is the simplest of the bilinear forms derived from the equation of continuity.
Therefore, it can be used as an indicator when considering several bilinear problems.
We will compare this equation with the model of fluid mechanics for example.
Usually, the nonlinearity of the bilinear convection-diffusion equation is given by $\boldsymbol{a}\cdot\nabla (|u|u)$ instead of $\boldsymbol{a}\cdot\nabla (u^2)$.
The reason for not choosing this form will be explained in the appendix.
Since \eqref{cd} comes from the equation of continuity, the conservation law $M_0 = \int_{\mathbb{R}^n} u_0 (x) dx = \int_{\mathbb{R}^n} u(t,x) dx$ holds.
Furthermore, from Moser--Nash method, we see that
\begin{equation}\label{decay}
	\| u(t) \|_{L^q (\mathbb{R}^n)} \le C(1+t)^{-\gamma_q}
\end{equation}
for $1 \le q \le \infty$ and $\gamma_q = \frac{n}2 (1-\frac1q)$ if $u_0 \in L^1 (\mathbb{R}^n) \cap L^\infty (\mathbb{R}^n)$.
Moreover, it is well known that there are unique functions $U_m$ such that
\[
	\lambda^{n+m} U_m (\lambda^2 t, \lambda x) = U_m (t,x)
\]
for $\lambda > 0$, and
\begin{equation}\label{EZexp}
	\biggl\| u(t) - \sum_{m=0}^{n-2} U_m (t) \biggr\|_{L^q (\mathbb{R}^n)}
	=
	O (t^{-\gamma_q-\frac{n}2+\frac12} \log t)
\end{equation}
as $t \to +\infty$.
More precisely, the linear components
\begin{equation}\label{Um}
\begin{split}
	&U_m (t) = \sum_{|\alpha|=m} \frac{\nabla^\alpha G(t)}{\alpha!} \int_{\mathbb{R}^n} (-y)^\alpha u_0 (y) dy\\
	&+ \sum_{2l+|\beta|=m-1} \frac{\partial_t^l \nabla^\beta (\boldsymbol{a}\cdot\nabla) G(t)}{l!\beta!} \int_0^\infty \int_{\mathbb{R}^n} (-s)^l (-y)^\beta u^2 (s,y) dyds
\end{split}
\end{equation}
for $0 \le m \le n-2$ are decided, where $G(t,x) = (4\pi t)^{-n/2} e^{-|x|^2/(4t)}$ is $n$-dimensional Gaussian.
For the details of these facts, see \cite{EZ,DrCrpo}.
Those terms represent the linear component of the solution.
Indeed, they fulfill that
\[
	\partial_t U_m = \Delta U_m.
\]
Therefore, they naturally obey the second law of thermodynamics.
In fact, the decay rate of $U_m$ is described by $m$, but simultanously $m$ also represents the symmetry of $U_m$.
What concerns us here is whether the logarithm appearing in \eqref{EZexp} is essential or not.
The more fundamental question is what is the largest nonlinear component.
To solve it, we will derive higher-order expansion, but first we need to confirm the convergence of the moments.
We see for the coefficients of the second part of $U_m$ that
\[
\begin{split}
	&\biggl| \int_0^\infty \int_{\mathbb{R}^n} (-s)^l (-y)^\beta u^2 (s,y) dyds \biggr|
	\le
	\int_0^\infty s^l \| y^\beta u (s) \|_{L^2 (\mathbb{R}^n)}^2 ds\\
	&\le
	C \int_0^\infty s^{\frac{m-1}2} (1+s)^{-n/2} ds,
\end{split}
\]
where we apply the weighted estimate that
$
	\| |x|^\mu u(t) \|_{L^q (\mathbb{R}^n)} = O (t^{-\gamma_q+\frac\mu2})
$ as $t \to +\infty$ (see the last sentence of this section).
The right-hand side diverges to infinity if $m \ge n-1$.
For this reason, the order of the expansion \eqref{EZexp} is limited.
The renormalization is affective to avoid this difficulty.
This method is well known as a technique for clarifying asymptotic structures of solutions to critical problems, but it can also be applied to derive higher-order asymptotic expansions of solutions to subcritical problems such as \eqref{cd}.
In the two-dimensional case of \eqref{cd}, the renormalization is synonymous with considering the linearized problem.
By studying the linearized problem, several researchers proved that the logarithmic evolution in \eqref{EZexp} is crucial.
\begin{theorem}[cf.\cite{EZ,DrCrpo,DrZZ}]\label{thm2dim}
Let $n = 2,~ u_0 \in L^1 (\mathbb{R}^2) \cap L^\infty (\mathbb{R}^2)$ and $|x| u_0 \in L^1 (\mathbb{R}^2)$.
Then, the solution $u$ of \eqref{cd} fulfills that
\[
	\left\| u(t) - M_0 G(t) - \frac{M_0^2}{8\pi} \boldsymbol{a}\cdot \nabla G(t) \log t \right\|_{L^q (\mathbb{R}^2)}
	= O(t^{-\gamma_q-\frac12})
\]
as $t \to +\infty$ for $1 \le q \le \infty$ and $\gamma_q = 1 - \frac1q$, where $M_0 = \int_{\mathbb{R}^2} u_0 (x) dx$.
\end{theorem}
In higher-dimensions, the renormalization is affective.
We deal with the four-dimensional case first, leaving the three-dimensional case for later.
In four dimensions, the estimate is still crucial.
\begin{theorem}\label{thm4dim}
Let $n = 4,~ u_0 \in L^1 (\mathbb{R}^4) \cap L^\infty (\mathbb{R}^4)$ and $|x|^3 u_0 \in L^1 (\mathbb{R}^4)$.
Then, the solution $u$ of \eqref{cd} fulfills that
\[
\begin{split}
	\biggl\| &u(t) - M_0 G(t) - \boldsymbol{M}_1\cdot\nabla G(t) - \sum_{|\alpha| = 2} \frac{\nabla^\alpha G(t)}{\alpha!} \int_{\mathbb{R}^4} (-y)^\alpha u_0 (y) dy\\
	&+ \nabla (\boldsymbol{a}\cdot\nabla) G(t) \cdot \int_0^\infty \int_{\mathbb{R}^4} y u^2 (s,y) dyds
	 + \frac{M_0^2}{128\pi^2} \Delta (\boldsymbol{a}\cdot\nabla) G(t) \log t \biggr\|_{L^q (\mathbb{R}^4)}
	= O(t^{-\gamma_q-\frac32})
\end{split}
\]
as $t \to +\infty$ for $1 \le q \le \infty$ and $\gamma_q = 2(1-\frac1q)$, where $M_0 = \int_{\mathbb{R}^4} u_0 (x) dx \in \mathbb{R}$ and $\boldsymbol{M}_1 = - \int_{\mathbb{R}^4} x u_0 (x) dx + a \int_0^\infty \int_{\mathbb{R}^4} u^2 (t,x) dxdt \in \mathbb{R}^4$.
\end{theorem}
The logarithmic evolutions in these theorems indicate the nonlinear distortions of symmetries in the solutions, respectively.
Indeed, when we denote these terms by $K_1 (t) = \frac{M_0^2}{8\pi} \boldsymbol{a}\cdot \nabla G(t)$ and $K_3(t) = - \frac{M_0^2}{128\pi^2} \Delta (\boldsymbol{a}\cdot\nabla) G(t)$, we see that $\lambda^{n+m} K_m (\lambda^2 t, \lambda x) = K_m (t,x)$ for $\lambda > 0$ and
\[
	\partial_t (K_m \log t) = \Delta (K_m \log t) - \frac{(K_m \log t)}{t\log t}.
\]
We also note that they are decided only by $a$ and $M_0 = \int_{\mathbb{R}^n} u_0 (x) dx$.
Based on \eqref{EZexp} and the above assertions, for the three-dimensional case, we can expect that there is some function $K_2$ such that $\lambda^{5} K_2 (\lambda^2 t, \lambda x) = K_2 (t,x)$ for $\lambda > 0$, and $u(t) \sim U_0 (t) + U_1 (t) + K_2 (t) \log t$ as $t \to +\infty$.
In fact, we find the distortion of other form in our main result as follows.
\begin{theorem}\label{thm3dim}
Let $n = 3,~ u_0 \in L^1 (\mathbb{R}^3) \cap L^\infty (\mathbb{R}^3)$ and $|x|^4 u_0 \in L^1 (\mathbb{R}^3)$.
Then, the solution $u$ of \eqref{cd} fulfills that
\[
\begin{split}
	\biggl\| &u(t) - M_0 G(t) - \boldsymbol{M}_1 \cdot \nabla G(t) - \sum_{|\alpha|=2}^3 \frac{\nabla^\alpha G(t)}{\alpha!} \int_{\mathbb{R}^3} (-y)^\alpha u_0 (y) dy\\
	&+ \nabla (\boldsymbol{a}\cdot\nabla) G(t) \cdot \int_0^\infty \int_{\mathbb{R}^3} y u^2 (s,y) dyds\\
	&- \sum_{2l+|\beta|=2} \frac{\partial_t^l \nabla^\beta (\boldsymbol{a}\cdot\nabla) G(t)}{\beta!} \int_0^\infty \int_{\mathbb{R}^3} (-s)^l (-y)^\beta (u^2 - M_0^2 G^2) (s,y) dyds\\
	&+ \frac{\sqrt{2\pi}}{32\pi^2} M_0 (M_0+\boldsymbol{M}_1\cdot\nabla) \left( 2t^{-1/2} \boldsymbol{a}\cdot\nabla G(\tfrac{t}2) + \int_0^t s^{-1/2} \Delta (\boldsymbol{a}\cdot\nabla) G(t-\tfrac{s}2) ds \right)\\
	&+ \frac{\sqrt{3} M_0^3}{2^6 \cdot 3^3 \cdot 5\pi^3} \Delta (\boldsymbol{a}\cdot\nabla)^2 G(t) \log t \biggr\|_{L^q (\mathbb{R}^3)}
	= O(t^{-\gamma_q-2})
\end{split}
\]
as $t \to +\infty$ for $1 \le q \le \infty$ and $\gamma_q = \frac32 (1-\frac1q)$, where $M_0 = \int_{\mathbb{R}^3} u_0 (x) dx \in \mathbb{R}$ and $\boldsymbol{M}_1 = - \int_{\mathbb{R}^3} x u_0 (x) dx + a \int_0^\infty \int_{\mathbb{R}^3} u^2 (t,x) dxdt \in \mathbb{R}^3$.
\end{theorem}
Namely, \eqref{EZexp} is not crucial in three dimensions, and there are unique functions $U_m$ and $K_4$ such that
\[
	\lambda^{3+m} U_m (\lambda^2 t, \lambda x) = U_m (t,x),\quad
	\lambda^{7} K_4 (\lambda^2 t, \lambda x) = K_4 (t,x)
\]
for $\lambda > 0$, and
\[
	u(t) \sim U_0 (t) + U_1 (t) + U_2 (t) + U_3 (t) + K_4 (t) \log t
\]
as $t \to +\infty$.
Here $U_0$ and $U_1$ are linear components defined in \eqref{Um}, i.e.,
\begin{equation}\label{U0U1}
	U_0 (t) = M_0 G(t),\quad
	U_1 (t) = \boldsymbol{M}_1 \cdot \nabla G(t).
\end{equation}
The higher-order terms are introduced here.
For example,
\begin{equation}\label{U2}
\begin{split}
	&U_2 (t) = \sum_{|\alpha|=2} \frac{\nabla^\alpha G(t)}{\alpha!} \int_{\mathbb{R}^3} y^\alpha u_0 (y) dy
	- \nabla (\boldsymbol{a}\cdot\nabla) G(t) \cdot \int_0^\infty \int_{\mathbb{R}^3} y u^2 (s,y) dyds\\
	&-  \frac{\sqrt{2\pi}}{32\pi^2} M_0^2 \left( 2t^{-1/2} \boldsymbol{a}\cdot\nabla G(\tfrac{t}2) + \int_0^t s^{-1/2} \Delta (\boldsymbol{a}\cdot\nabla) G(t-\tfrac{s}2) ds \right).
\end{split}
\end{equation}
In this term, we find the nonlinear distortion in the solution symmetry since $\partial_t U_2 = \Delta U_2$ is never satisfied.
More precisely, the former two and last parts of $U_2$ have different symmetries.
Namely, if we put $U_2 = U_2^{\mathrm{evn}} + U_2^{\mathrm{odd}}$, then $U_2^{\mathrm{evn}} (t,-x) = U_2^{\mathrm{evn}} (t,x)$ and $U_2^{\mathrm{odd}}(t,-x) = - U_2^{\mathrm{odd}} (t,x)$ hold for $(t,x) \in \mathbb{R}_+ \times \mathbb{R}^3$.
This structure plays important role in the derivation process of the logarithmic term.
In linear problems, due to the second law of thermodynamics, such parities corresponds to the parabolic scale.
Indeed, the linear components $U_0$ and $U_1$ are even and odd-type functions, respectively.
Specifically, $U_1$ exhibits higher symmetry than $U_0$.
Conversely, $U_2$ includes several parities.
To be clear, $U_2^{\mathrm{evn}}$ is a linear component essentially since $\partial_t U_2^{\mathrm{evn}} = \Delta U_2^{\mathrm{evn}}$ still holds, and $U_2^{\mathrm{odd}}$ contains both component with higher symmetry and component with lower symmetry compared to $U_2^{\mathrm{evn}}$.
The first half $2 t^{-1/2} \boldsymbol{a}\cdot\nabla G(\tfrac{t}2)$ decays rapidly despite its low symmetry, while the latter half $\int_0^t s^{-1/2} \Delta (\boldsymbol{a}\cdot\nabla) G(t-\tfrac{s}2) ds$ decays slowly despite its high symmetry.
This is the largest profile where linear and nonlinear effects are balanced.
The next term
\begin{equation}\label{U3}
\begin{split}
	&U_3 (t) = -\sum_{|\alpha| = 3} \frac{\nabla^\alpha G(t)}{\alpha!} \int_{\mathbb{R}^3} y^\alpha u_0 (y) dy\\
	&+ \sum_{2l+|\beta|=2} \frac{\partial_t^l \nabla^\beta (\boldsymbol{a}\cdot\nabla) G(t)}{\beta!} \int_0^\infty \int_{\mathbb{R}^3} (-s)^l (-y)^\beta (u^2 - M_0^2 G^2) (s,y) dyds\\
	&- \frac{\sqrt{2\pi}}{32\pi^2} M_0 \left( 2t^{-1/2} (\boldsymbol{M}_1\cdot\nabla) (\boldsymbol{a}\cdot\nabla) G(\tfrac{t}2) + \int_0^t s^{-1/2} \Delta (\boldsymbol{M}_1\cdot\nabla) (\boldsymbol{a}\cdot\nabla) G (t-\tfrac{s}2) ds \right)
\end{split}
\end{equation}
has the similar structure.
The logarithmic evolution is appearing as a tiny distortion but its coefficient
\begin{equation}\label{K4}
	K_4 (t) = -\frac{\sqrt{3} M_0^3}{2^6 \cdot 3^3 \cdot 5\pi^3} \Delta (\boldsymbol{a}\cdot\nabla)^2 G(t)
\end{equation}
is also decided only by $a$ and $M_0 = \int_{\mathbb{R}^3} u_0 (x) dx$.
Reader may doubt that some moments in the expansion do not converge.
For instance, the coefficient of $U_2$ is estimated as
\[
	\biggl| \int_0^\infty \int_{\mathbb{R}^3} y u^2 (s,y) dyds \biggr|
	\le
	C \int_0^\infty s^{-1/2} (1+s)^{-1/2} ds = +\infty
\]
at first glance.
In fact, the renormalization fills this gap.
Specifically, this moment is rewritten by $\int_{\mathbb{R}^3} yu^2 dy = \int_{\mathbb{R}^3} y(u^2 - M_0^2 G^2) dy$ since $G^2$ is radially symmetric in $y$.
Hence, from \eqref{EZexp} and the weighted estimate, we have that
\[
	\biggl| \int_0^\infty \int_{\mathbb{R}^3} y u^2 (s,y) dyds \biggr|
	\le
	C \int_0^\infty s^{-1/2} (1+s)^{-1} ds < +\infty.
\]
The other moments are treated on the same way.

Also in the case $n=2$ and $n=4$, the solutions contain the similar distortions as in Theorem \ref{thm3dim}.
However, due to the influence of large logarithms, the concrete forms of them are complicated.
As are the case of two and four dimensions,  in general even-dimensional cases, \eqref{EZexp} is expected to be crucial and there is some unique function $K_{n-1}$ such that $\lambda^{2n-1} K_{n-1} (\lambda^2 t, \lambda x) = K_{n-1} (t,x)$ for $\lambda > 0$ and $u(t) \sim U_0 (t) + \cdots + U_{n-2} (t) + K_{n-1} (t) \log t$ as $t \to +\infty$ and the last term is the distortion.
On the other hand, in the odd-dimensional cases, the following expansion is derived.
\begin{proposition}\label{propoddim}
Let $n \ge 3$ be odd, $u_0 \in L^1 (\mathbb{R}^n) \cap L^\infty (\mathbb{R}^n)$ and $|x|^{2n-2} u_0 \in L^1 (\mathbb{R}^n)$.
Then there are unique functions $U_m$ such that $\lambda^{n+m} U_m (\lambda^2 t,\lambda x) = U_m (t,x)$ for $\lambda > 0$, and
\[
	\biggl\| u(t) - \sum_{m=0}^{2n-3} U_m (t) \biggr\|_{L^q (\mathbb{R}^n)}
	= O (t^{-\gamma_q-n+1} \log t)
\]
as $t \to +\infty$ for $1 \le q \le \infty$ and $\gamma_q = \frac{n}2 (1-\frac1q)$, and the solution $u$ of \eqref{cd}.
\end{proposition}
Of coursely, the lower-orders $U_0, U_1,\ldots,U_{n-2}$ are given by \eqref{Um}.
We will decide the higher-orders $U_{n-1},U_n,\ldots,U_{2n-3}$ of distortions.
In the three dimensions, the better conclusion is given in Theorem \ref{thm3dim}.
Even in the case $n \ge 5$, the solution may have some $K_{2n-2} (t) \log t$ as a profile.
In other words, the logarithmic components in the case of odd dimensions decay much faster than expected.
We should confirm whether $K_{2n-2}$ is remaining or vanishing.
Unfortunately, it is difficult to prove this for higher dimensions at present.
Note that the one-dimensional Burgers equation is scale-critical, and the behavior of its solutions differs significantly from that in other dimensions.
Indeed, it is proved by Kato \cite{KtM} that $u(t) \sim U_0(t) + K_1(t) \log t$ as $t \to +\infty$.
Here $U_0$ is not Gaussian but Cole--Hopf solution.
In this critical case, the both of $U_0$ and $K_1 \log t$ describe the nonlinear distortion.

Generally, the nonlinear term of convection-diffusion equation is given by $\boldsymbol{a}\cdot \nabla (|u|^{r-1}u)$ for some $r$.
Especially, the case $r = 1 + \frac1n$ is known as critical.
Both of asymptotic profile for the critical case and asymptotic expansion of the subcritical case are studied by several authors (cf.\cite{BKL,BKW,FkdSt,IshgKwkm13,Ksb}).
The word `renormalization' is often seen in the critical case.
Our equation corresponds to the case $r = 2$ and this case is subcritical.

The reason why the Burgers type \eqref{cd} is focused is that this equation is an indicator for considering other bilinear problems.
We refer to the incompressible Navier--Stokes flow for example.
As well known, this flow is formulated as
\begin{equation}\label{ns}
	\left\{\begin{array}{lr}
		\partial_t \boldsymbol{v} + \boldsymbol{v}\cdot\nabla \boldsymbol{v} = \Delta \boldsymbol{v} - \nabla p,	&	t>0,~ x \in \mathbb{R}^n,\\
		\nabla\cdot \boldsymbol{v} = 0,															&	t>0,~ x \in \mathbb{R}^n,\\
		\boldsymbol{v}(0,x) = \boldsymbol{v}_0 (x),												&	x \in \mathbb{R}^n,
	\end{array}\right.
\end{equation}
where the velocity $\boldsymbol{v}$, unlike $u$, is a vector.
The effects of pressure $p$ can be eliminated by the proper projection.
The solenoidal condition $\nabla\cdot \boldsymbol{v}_0 = 0$ is comparable to the loss of mass $\int_{\mathbb{R}^n} u_0 (x) dx = 0$ in our problem.
This leads to the expectation that the nonlinearity of the Navier--Stokes equation is two steps weaker than that of our equation.
Indeed, Fujigaki and Miyakawa\cite{FjgkMykw} derived the unique linear profiles $\boldsymbol{V}_m$ such that $\lambda^{n+m} \boldsymbol{V}_m (\lambda^2 t, \lambda x) = \boldsymbol{V}_m (t,x)$ for $\lambda > 0$, and
\[
	\biggl\| \boldsymbol{v}(t) - \sum_{m=1}^n \boldsymbol{V}_m (t) \biggr\|_{L^q (\mathbb{R}^n)}
	=
	O (t^{-\gamma_q-\frac{n}2 - \frac12} \log t)
\]
as $t \to +\infty$.
Here $\boldsymbol{V}_0$ is erased by the solenoidal condition.
If we restrict ourselves to even dimensions, a comparison of this estimate and \eqref{EZexp} shows that the above prediction is correct.
Theorems \ref{thm2dim} and \ref{thm4dim} support this fact.
Recently, for the odd dimensions, the author showed that
\[
	\biggl\| \boldsymbol{v}(t) - \sum_{m=1}^{2n} \boldsymbol{V}_m (t) \biggr\|_{L^q (\mathbb{R}^n)}
	=
	O (t^{-\gamma_q- n - \frac12} \log t)
\]
as $t \to +\infty$ (cf.\cite{Ym25}).
Here $\boldsymbol{V}_m$ for $1 \le m \le n$ are the linear profiles and one for $n+1 \le m \le 2n$ contains the nonlinear distortion.
Comparing this estimate with Proposition \ref{propoddim}, in the odd dimensions, we see that the logarithmic evolution in Navier--Stokes flow is at least three steps smaller than one in our solution.
Moreover, Theorem \ref{thm3dim} suggests that this difference is essential.
The reason why such a difference appears is that, in the odd-dimensional Navier--Stokes flows, the solenoidal condition changes not only the scale but also the symmetry of nonlinearity.
More closely related to our motivation is the following problem proposed by Carpio:
\begin{itemize}
\item
If fluid outflow or inflow occurs at the initial time, how does it affect the flow velocity?
\end{itemize}
This corresponds to considering the asymptotic expansion of $\boldsymbol{v}$ when $\nabla \cdot \boldsymbol{v}_0 \neq 0$.
Under this condition, \eqref{ns} has same scale as \eqref{cd}.
The suction and emission of fluid correspond to the linear waves in our equations.
In \cite{Crpo}, this problem of  two dimensions is studied and the similar assertion as Theorem \ref{thm2dim} is derived.
Namely, the effect of the fluid outflow or inflow distort the velocity and this distiortion evolvs logarithmically in time.
For the higher dimensions, due to the complexity of the symmetry, this problem remains unsolved.
As a clue to solving it, we introduced \eqref{cd} which has same scale as \eqref{ns} and possesses the simplest symmetry.
Based on Theorem \ref{thm3dim}, distortion of three-dimensional flow velocity is expected to be quite complex.

Before closing this section, we rewrite our problem and confirm the weighted estimate.
Additionally, we explain the asymptotic analysis of the linear problem.
Duhamel principle leads the mild solution of \eqref{cd} that
\begin{equation}\label{ms}
	u(t) = G(t) * u_0 + \int_0^t \boldsymbol{a}\cdot\nabla G(t-s) * u^2 (s) ds,
\end{equation}
where $*$ means spatial convolution and we omit the spatial variable.
For some $\mu \ge 0$, the term of initial-data fulfills that $\| |x|^\mu G(t)*u_0 \|_{L^q (\mathbb{R}^n)} \le Ct^{-\gamma_q} (1+t)^{\mu/2}$ for $1 \le q \le \infty$ if $|x|^\mu u_0 \in L^1 (\mathbb{R}^n)$ is assumed.
We put $f_q (t) = \sup_{0 < s < t} s^{\gamma_q} (1+s)^{-\mu/2} \| |x|^\mu u(s) \|_{L^q (\mathbb{R}^n)}$, then we have from \eqref{ms} with Hausdorff--Young inequality and \eqref{decay} that
\[
\begin{split}
	&\| |x|^\mu u(t) \|_{L^q (\mathbb{R}^n)}
	\le
	\| |x|^\mu G(t)*u_0 \|_{L^q (\mathbb{R}^n)}\\
	&+ C \int_0^t \| |x|^\mu (\boldsymbol{a}\cdot\nabla) G(t-s) \|_{L^q (\mathbb{R}^n)} \| u(s) \|_{L^2 (\mathbb{R}^n)}^2 ds\\
	&+ C \int_0^{t/2} \| (\boldsymbol{a}\cdot\nabla) G(t-s) \|_{L^q (\mathbb{R}^n)} \| |x|^\mu u^2 (s) \|_{L^2 (\mathbb{R}^n)} ds\\
	&+ C \int_{t/2}^t \| (\boldsymbol{a}\cdot\nabla) G(t-s) \|_{L^1 (\mathbb{R}^n)} \| |x|^\mu u^2 (s) \|_{L^q (\mathbb{R}^n)} ds\\
	&\le
	Ct^{-\gamma_q} (1+t)^{\mu/2} + C \int_0^t (t-s)^{-\gamma_q-\frac12+\frac\mu2} (1+s)^{-\frac{n}2} ds\\
	&+ C f_1 (t) \int_0^{t/2} (t-s)^{-\gamma_q-\frac12} (1+s)^{-\frac{n}2+\frac\mu2} ds\\
	&+ Cf_q (t) \int_{t/2}^t (t-s)^{-\frac12} s^{-\gamma_q} (1+s)^{-\frac{n}2+\frac\mu2} ds.
\end{split}
\]
This leads that
\begin{equation}\label{decay-wt}
	\| |x|^\mu u(t) \|_{L^q (\mathbb{R}^n)} \le Ct^{-\gamma_q} (1+t)^{\mu/2}
\end{equation}
for $1 \le q \le \infty$ when $|x|^\mu u_0 \in L^1 (\mathbb{R}^n)$ is supposed.
The asymptotic expansion is derived from \eqref{ms} by applying Taylor theorem to its integral kernels.
For example, the term of initial-data is expanded as
\[
	G(t) * u_0 = \sum_{|\alpha|=0}^{m-1} \frac{\nabla^\alpha G(t)}{\alpha!} \int_{\mathbb{R}^n} (-y)^\alpha u_0 (y) dy + r_m^0 (t)
\]
for
\begin{equation}\label{rm0}
	r_m^0 (t) = \int_{\mathbb{R}^n} \biggl( G(t,x-y) - \sum_{|\alpha|=0}^{m-1} \frac{\nabla^\alpha G(t,x)}{\alpha!} (-y)^\alpha \biggr) u_0 (y) dy.
\end{equation}
Any $U_m$ in our assertions contain this form and Taylor theorem says that $\| r_m^0 (t) \|_{L^q (\mathbb{R}^n)} = O (t^{-\gamma_q-\frac{m}2})$ as $t \to +\infty$ for $1 \le q \le \infty$ if $|x|^m u_0 \in L^1 (\mathbb{R}^n)$ is assumed.

\vspace{2mm}

\paragraph{\textbf{Notations.}}
We often omit the spatial parameter from functions, for example, $u(t) = u(t,x)$.
In particular, $G(t) * u_0 = \int_{\mathbb{R}^n} G(t,x-y) u_0 (y) dy$ and $\int_0^t g(t-s) * f(s) ds = \int_0^t \int_{\mathbb{R}^n} g(t-s,x-y) f(s,y) dyds$.
We symbolize the derivations by $\partial_t = \partial/\partial t,~ \partial_j = \partial/\partial x_j$ for $1 \le j \le n,~ \nabla = (\partial_1,\partial_2,\ldots,\partial_n)$ and $\Delta = \lvert \nabla \rvert^2 = \partial_1^2 + \partial_2^2 + \cdots + \partial_n^2$.
The length of a multiindex $\alpha = (\alpha_1,\alpha_2,\ldots, \alpha_n) \in \mathbb{Z}_+^n$ is given by $\lvert \alpha \rvert = \alpha_1 + \alpha_2 + \cdots + \alpha_n$, where $\mathbb{Z}_+ = \mathbb{N} \cup \{ 0 \}$.
We abbreviate that $\alpha ! = \alpha_1 ! \alpha_2! \cdots \alpha_n !,~ x^\alpha = x_1^{\alpha_1} x_2^{\alpha_2} \cdots x_n^{\alpha_n}$ and $\nabla^\alpha = \partial_1^{\alpha_1} \partial_2^{\alpha_2} \cdots \partial_n^{\alpha_n}$.
We define the Fourier transform and its inverse by $\mathcal{F} [\varphi] (\xi) = (2\pi)^{-n/2} \int_{\mathbb{R}^n} \varphi (x) e^{-ix\cdot\xi} dx$ and $\mathcal{F}^{-1} [\varphi] (x) = (2\pi)^{-n/2} \int_{\mathbb{R}^n} \varphi (\xi) e^{ix\cdot\xi} d\xi$, respectively, where $i = \sqrt{-1}$.
The Lebesgue space and its norm are denoted by $L^q (\mathbb{R}^n)$ and $\| \cdot \|_{L^q (\mathbb{R}^n)}$, that is, $\| f \|_{L^q (\mathbb{R}^n)} = (\int_{\mathbb{R}^n} |f(x)|^q dx)^{1/q}$ for $1 \le q < \infty$ and $\| f \|_{L^\infty (\mathbb{R}^n)}$ is the essential supremum.
The heat kernel and its decay rate on $L^q (\mathbb{R}^n)$ are symbolized by $G(t,x) = (4\pi t)^{-n/2} e^{-|x|^2/(4t)}$ and $\gamma_q = \frac{n}2 (1-\frac1q)$.
We employ Landau symbol.
Namely, $f(t) = o(t^{-\mu})$ and $g(t) = O(t^{-\mu})$ mean $t^\mu f(t) \to 0$ and $t^\mu g(t) \to c$ for some $c \in \mathbb{R}$ such as $t \to +\infty$ or $t \to +0$, respectively.
A subscript of function represents its scale or decay rate.
For example, $\| |x|^\mu U_m (t) \|_{L^q (\mathbb{R}^n)} = t^{-\gamma_q-\frac{m}2+\frac\mu2} \| |x|^\mu U_m (1) \|_{L^q (\mathbb{R}^n)}$ for $t > 0$, and $\| r_m (t) \|_{L^q (\mathbb{R}^n)} = O (t^{-\gamma_q-\frac{m}2})$ or $O(t^{-\gamma_q-\frac{m}2} \log t)$ as $t \to +\infty$.
Various positive constants are simply denoted by $C$.

\section{The even-dimensional cases}
We prove Theorems \ref{thm2dim} and \ref{thm4dim}.
\subsection{The two-dimensional case}
This case is well-known and reader may skip it.
We expand the right-hand side of \eqref{ms} as
\begin{equation}\label{bs2dim}
\begin{split}
	&u(t) = G(t) \int_{\mathbb{R}^2} u_0 (y) dy
	+ \boldsymbol{a}\cdot\nabla G(t) \int_0^t \int_{\mathbb{R}^2} u^2 (s,y) dyds
	+ r_1^0 (t) +  r_1^1 (t)
\end{split}
\end{equation}
for $r_1^0$ given by \eqref{rm0}
and
\[
	r_1^1 (t) = \int_0^t \int_{\mathbb{R}^2} \left( \boldsymbol{a}\cdot\nabla G(t-s,x-y) - \boldsymbol{a}\cdot\nabla G(t,x) \right) u^2 (s,y) dyds.
\]
The first term is $U_0$ and the other terms decay fast as $t \to +\infty$.
The mean value theorem leads that
\[
\begin{split}
	&r_1^1 (t) = - \int_0^{t/2} \int_0^1 \partial_t (\boldsymbol{a}\cdot\nabla) G(t-\lambda s) * s u^2 (s) d\lambda ds\\
	&- \int_0^{t/2} \int_{\mathbb{R}^2} \int_0^1 (y\cdot\nabla) (\boldsymbol{a}\cdot\nabla) G(t, x-\lambda y) u^2 (s,y) d\lambda dyds\\
	&+
	\int_{t/2}^t \int_{\mathbb{R}^2} (\boldsymbol{a}\cdot\nabla G(t-s,x-y) - \boldsymbol{a}\cdot\nabla G(t,x)) u^2 (s,y) dyds.
\end{split}
\]
For instance, from \eqref{decay} and \eqref{decay-wt},
\[
\begin{split}
	&\biggl\| \int_0^{t/2} \int_0^1 \partial_t (\boldsymbol{a}\cdot\nabla) G(t-\lambda s) * s u^2 (s) d\lambda ds \biggr\|_{L^q (\mathbb{R}^2)}\\
	&+\biggl\| \int_0^{t/2} \int_{\mathbb{R}^2} \int_0^1 (y\cdot\nabla) (\boldsymbol{a}\cdot\nabla) G(t, x-\lambda y) u^2 (s,y) d\lambda dyds \biggr\|_{L^q (\mathbb{R}^2)}\\
	&\le
	\int_0^{t/2} \int_0^1 \| \partial_t (\boldsymbol{a}\cdot\nabla) G(t-\lambda s) \|_{L^q (\mathbb{R}^2)} s^l \| u^2 (s) \|_{L^1 (\mathbb{R}^2)} d\lambda ds\\
	&+
	\int_0^{t/2} \int_0^1 \| \nabla (\boldsymbol{a}\cdot\nabla) G(t) \|_{L^q (\mathbb{R}^2)} \| xu^2 (s) \|_{L^1 (\mathbb{R}^2)} d\lambda ds
	\le
	C t^{-\gamma_q-\frac12}.
\end{split}
\]
Another part is treated as
\[
\begin{split}
	&\biggl\| \int_{t/2}^t (\boldsymbol{a}\cdot\nabla G(t-s,x-y) - \boldsymbol{a}\cdot\nabla G(t,x)) u^2 (s,y) dyds \biggr\|_{L^q (\mathbb{R}^n)}\\
	&\le
	\int_{t/2}^t \| \boldsymbol{a}\cdot\nabla G(t-s) \|_{L^1 (\mathbb{R}^n)} \| u^2(s) \|_{L^q (\mathbb{R}^n)} ds
	+
	\int_{t/2}^t \| \boldsymbol{a}\cdot\nabla G(t) \|_{L^q (\mathbb{R}^n)} \| u(s) \|_{L^2 (\mathbb{R}^n)}^2 ds\\
	&\le
	C t^{-\gamma_q-\frac12}.
\end{split}
\]
Namely, $\| r_1^1 (t) \|_{L^q (\mathbb{R}^3)} = O (t^{-\gamma_q-\frac12})$ as $t \to +\infty$.
The logarithmic evolution comes from the second term of \eqref{bs2dim}.
Indeed, by renormalizing $u^2$ by $M_0^2 G^2$, we see
\[
\begin{split}
	&\int_0^t \int_{\mathbb{R}^2} u^2 (s,y) dyds
	=
	M_0^2 \int_0^t \int_{\mathbb{R}^2} G^2 (1+s,y) dyds\\
	&+
	\int_0^t \int_{\mathbb{R}^2} (u^2 (s,y) - M_0^2 G^2(1+s,y)) dyds.
\end{split}
\]
The second term is uniformly integrable from \eqref{EZexp}, and the parabolic scale of $G$ yields that
\[
	\int_0^t \int_{\mathbb{R}^2} G^2 (1+s,y) dyds
	=
	\int_0^t (1+s)^{-1} ds \int_{\mathbb{R}^2} G^2 (1,y) dy
	=
	\frac1{8\pi} \log (1+t).
\]
Therefore we complete the proof.
\subsection{The four-dimensional case}\label{subsec4dim}
The proof is based on the same procedure as above with Taylor theorem instead of the mean value theorem.
Namely, we see that
\[
\begin{split}
	&u(t) = \sum_{|\alpha| = 0}^2 \frac{\nabla^\alpha G(t)}{\alpha!} \int_{\mathbb{R}^4} (-y)^\alpha u_0 (y) dy\\
	&+ \sum_{2l+|\beta|=0}^2 \frac{\partial_t^l \nabla^\beta (\boldsymbol{a}\cdot\nabla) G(t)}{\beta!} \int_0^t \int_{\mathbb{R}^4} (-s)^l (-y)^\beta u^2 (s,y) dyds
	+ r_3^0 (t) + r_3^1 (t)
\end{split}
\]
for $r_3^0$ given by \eqref{rm0}
and
\[
\begin{split}
	&r_3^1 (t)	= \int_0^t \int_{\mathbb{R}^4} \biggl( \boldsymbol{a}\cdot\nabla G(t-s,x-y) - \sum_{2l+|\beta|=0}^2 \frac{\partial_t^l \nabla^\beta G(t,x)}{\beta!} (-s)^l (-y)^\beta \biggr) u^2 (s,y) dyds.
\end{split}
\]
For $|\beta| \le 1$, the coefficent of the second part is uniformly integrable in time.
Hence, we see
\[
	\int_0^t \int_{\mathbb{R}^4} (-y)^\beta u^2 (s,y) dyds
	=
	\int_0^\infty \int_{\mathbb{R}^4} (-y)^\beta u^2 (s,y) dyds - \int_t^\infty \int_{\mathbb{R}^4} (-y)^\beta u^2 (s,y) dyds.
\]
For $2l+|\beta|=2$, we separate the coefficient as
\[
\begin{split}
	&\int_0^t \int_{\mathbb{R}^4} (-s)^l (-y)^\beta u^2 (s,y) dyds
	= M_0^2 \int_0^t \int_{\mathbb{R}^4} (-s)^l (-y)^\beta G^2 (1+s,y) dyds\\
	&+ \int_0^t \int_{\mathbb{R}^4} (-s)^l (-y)^\beta (u^2 (s,y) - M_0^2 G^2 (1+s,y)) dyds.
\end{split}
\]
Finally, we have that
\[
\begin{split}
	&u(t) = \sum_{|\alpha| = 0}^2 \frac{\nabla^\alpha G(t)}{\alpha!} \int_{\mathbb{R}^4} (-y)^\alpha u_0 (y) dy\\
	&+ \sum_{|\beta|=0}^1 \frac{\nabla^\beta (\boldsymbol{a}\cdot\nabla) G(t)}{\beta!} \int_0^\infty \int_{\mathbb{R}^4} (-y)^\beta u^2 (s,y) dyds\\
	&+ M_0^2 \sum_{2l+|\beta|=2} \frac{\partial_t^l \nabla^\beta (\boldsymbol{a}\cdot\nabla) G(t)}{\beta!} \int_0^t \int_{\mathbb{R}^4} (-s)^l (-y)^\beta G^2 (1+s,y) dyds\\
	&+ r_3^0 (t) + r_3^1 (t) + r_3^2 (t) + r_3^3 (t)
\end{split}
\]
for
\[
	r_3^2 (t) = -  \sum_{|\beta|=0}^1 \frac{\nabla^\beta (\boldsymbol{a}\cdot\nabla) G(t)}{\beta!} \int_t^\infty \int_{\mathbb{R}^4} (-y)^\beta u^2 (s,y) dyds
\]
and
\[
	r_3^3 (t) = \sum_{2l+|\beta|=2} \frac{\partial_t^l \nabla^\beta (\boldsymbol{a}\cdot\nabla) G(t)}{\beta!} \int_0^t \int_{\mathbb{R}^4} (-s)^l (-y)^\beta \left( u^2 (s,y) - M_0^2 G^2 (1+s,y) \right) dyds.
\]
The first and second parts provide $U_0, U_1$ and $U_2$.
The remaind terms $r_3^0$ and $r_3^1$ could be estimated on the same way as above and we see $\| r_3^0 (t) \|_{L^q (\mathbb{R}^4)} + \| r_3^1 (t) \|_{L^q (\mathbb{R}^4)} = O (t^{-\gamma_q-\frac32})$ as $t\to\infty$.
The coefficients of $r_3^2$ fulfill that
\[
\begin{split}
	&\biggl| \int_t^\infty \int_{\mathbb{R}^4} (-y)^\beta u^2 (s,y) dyds \biggr|
	\le
	\int_t^\infty \| y^\beta u^2 (s) \|_{L^1 (\mathbb{R}^4)} ds\\
	&\le
	C \int_t^\infty s^{-2+\frac{|\beta|}2} ds
	= Ct^{-1+\frac{|\beta|}2}.
\end{split}
\]
Thus, $r_3^2$ also is an error term.
On the other hand, from \eqref{EZexp} and \eqref{decay-wt}, the coefficient of $r_3^3$ is uniformly integrable.
Thus, $\| r_3^2 (t) \|_{L^q (\mathbb{R}^4)} + \| r_3^3 (t) \|_{L^q (\mathbb{R}^4)} = O (t^{-\gamma_q-\frac32})$ as $t\to\infty$.
The logarithmic evolution is hyden in the third part.
Indeed, the coefficient is further separated as
\[
	\int_0^t \int_{\mathbb{R}^4} (-s)^l (-y)^\beta G^2 (1+s,y) dyds
	=
	\int_0^t s^l (1+s)^{-l-1} ds \int_{\mathbb{R}^4} (-1)^l (-y)^\beta G^2 (1,y) dy
\]
since $|\beta| = 2 - 2l$.
The integration in time yields the logarithmic evolution, that is $\int_0^t s^l (1+s)^{-l-1} ds = \log t + O(1)$ as $t \to\infty$, but the corresponding spatial integration is vanishing for $\beta_j = \beta_k = 1$ for some $j \neq k$.
Therefore, the coefficient of the logarithmic evolution is given by
\[
\begin{split}
	&M_0^2 \sum_{2l+|\beta|=2} \frac{\partial_t^l \nabla^\beta (\boldsymbol{a}\cdot\nabla) G(t)}{\beta!} \int_{\mathbb{R}^4} (-1)^l (-y)^\beta G^2 (1,y) dy\\
	&=
	M_0^2 \sum_{j=1}^4 \frac{\partial_j^2 (\boldsymbol{a}\cdot\nabla) G(t)}{2!} \int_{\mathbb{R}^4} y_j^2 G^2 (1,y) dy
	 - M_0^2\partial_t (\boldsymbol{a}\cdot\nabla) G(t) \int_{\mathbb{R}^4} G^2 (1,y) dy\\
	&= -\frac{M_0^2}{128\pi^2} \Delta (\boldsymbol{a}\cdot\nabla) G(t).
\end{split}
\]
Thus, we complete the proof.
\section{The three-dimensional case}\label{sect3dim}
We prove Theorem \ref{thm3dim} in this section.
Since this process is little complicated, we separate it to three steps.
\subsection{The expansion up to second order}\label{subsect3dim1}
Firstly, we expand $u$ up to second order.
By the similar argument as in Section \ref{subsec4dim}, we see that
\begin{equation}\label{bs1}
\begin{split}
	&u(t) = \sum_{|\alpha| = 0}^2 \frac{\nabla^\alpha G(t)}{\alpha!} \int_{\mathbb{R}^3} (-y)^\alpha u_0 (y) dy\\
	&+ \sum_{|\beta|=0}^1 \nabla^\beta (\boldsymbol{a}\cdot\nabla) G(t) \int_0^\infty \int_{\mathbb{R}^3} (-y)^\beta u^2 (s,y) dyds\\
	&+ \int_0^t \int_{\mathbb{R}^3} \biggl( \boldsymbol{a}\cdot\nabla G(t-s,x-y) - \sum_{|\beta|=0}^1 \nabla^\beta (\boldsymbol{a}\cdot\nabla) G(t,x) (-y)^\beta \biggr) u^2 (s,y) dyds\\
	&- \sum_{|\beta|=0}^1 \nabla^\beta (\boldsymbol{a}\cdot\nabla) G(t) \int_t^\infty \int_{\mathbb{R}^3} (-y)^\beta u^2 (s,y) dyds	
	+ r_3^0 (t)
\end{split}
\end{equation}
for $r_3^0$ given by \eqref{rm0}.
The first and second terms are $U_0, U_1$ written by \eqref{U0U1} and the part of $U_2$ introduced as \eqref{U2}.
Here, some components diverge to infinity at first glance.
In fact, the coefficients of the second part are already treated after the statement of this theorem.
On the similar way, we will treat the fourth part later.
We clarify large-time behavior of the third part of \eqref{bs1}.
By renormalizing $u^2$ by $M_0^2 G^2 = U_0^2$, we derive that
\begin{equation}\label{bs13}
\begin{split}
	&\int_0^t \int_{\mathbb{R}^3} \biggl( \boldsymbol{a}\cdot\nabla G(t-s,x-y) - \sum_{|\beta|=0}^1 \nabla^\beta (\boldsymbol{a}\cdot\nabla) G(t,x) (-y)^\beta \biggr) u^2 (s,y) dyds\\
	&=
	J_2 (t) + r_3^1 (t)
\end{split}
\end{equation}
for
\[
	J_2 (t) = M_0^2 \int_0^t \int_{\mathbb{R}^3} \biggl( \boldsymbol{a}\cdot\nabla G(t-s,x-y) - \sum_{|\beta|=0}^1 \nabla^\beta (\boldsymbol{a}\cdot\nabla) G(t,x) (-y)^\beta \biggr) G^2 (s,y) dyds
\]
and
\begin{equation}\label{r31}
	r_3^1 (t) = \int_0^t \int_{\mathbb{R}^3} \biggl( \boldsymbol{a}\cdot\nabla G(t-s,x-y) - \sum_{|\beta|=0}^1 \nabla^\beta (\boldsymbol{a}\cdot\nabla) G(t,x) (-y)^\beta \biggr) (u^2-M_0^2G^2) (s,y) dyds.
\end{equation}
%
Remark that this $r_3^1$ is similar as one in Section \ref{subsec4dim} but, since the parities are different, their definitions are little different.
One may doubt that $G^2 (s,y)$ in $J_2$ has a singularity as $s \to +0$.
In fact, by Taylor theorem, we obtain that
\[
\begin{split}
	&\int_0^{t/2} \int_{\mathbb{R}^3} \biggl( \boldsymbol{a}\cdot\nabla G(t-s,x-y) - \sum_{|\beta|=0}^1 \nabla^\beta (\boldsymbol{a}\cdot\nabla) G(t,x) (-y)^\beta \biggr) G^2 (s,y) dyds\\
	&=
	-\int_0^{t/2} \int_0^1 \partial_t (\boldsymbol{a}\cdot\nabla) G(t-\lambda s) s  G^2 (s) d\lambda ds\\
	&+
	\sum_{|\beta|=2} \int_0^{t/2} \int_{\mathbb{R}^3} \int_0^1 \frac{\nabla^\beta (\boldsymbol{a}\cdot\nabla) G(t, x-\lambda y)}{\beta!} \lambda (-y)^\beta G^2 (s,y) d\lambda dyds.
\end{split}
\]
Hence, Lebesgue convergence theorem solves the singularity.
From Taylor theorem, we see that
\[
\begin{split}
	&r_3^1 (t)
	= -\int_0^{t/2} \int_{\mathbb{R}^3} \int_0^1 \partial_t (\boldsymbol{a}\cdot\nabla) G(t-\lambda s) *s (u^2 - M_0^2G^2) (s) d\lambda ds\\
	&+
	\sum_{|\beta|=2} \int_0^{t/2} \int_{\mathbb{R}^3} \int_0^1 \frac{\nabla^\beta (\boldsymbol{a}\cdot\nabla) G(t, x-\lambda y)}{\beta!} \lambda y^\beta (u^2 - M_0^2G^2) (s,y) d\lambda dyds\\
	&+
	\int_{t/2}^t \int_{\mathbb{R}^3} \biggl( \boldsymbol{a}\cdot\nabla G(t-s,x-y) - \sum_{|\beta|=0}^1 \nabla^\beta (\boldsymbol{a}\cdot\nabla) G(t,x) (-y)^\beta \biggr)\\
	&\hspace{20mm} (u^2-M_0^2G^2) (s,y) dyds,
\end{split}
\]
and Hausdorff--Young inequality leads
\[
\begin{split}
	&\biggl\| \int_0^{t/2} \int_{\mathbb{R}^3} \int_0^1 \partial_t (\boldsymbol{a}\cdot\nabla) G(t-\lambda s) *s (u^2 - M_0^2G^2) (s) d\lambda ds \biggr\|_{L^q (\mathbb{R}^3)}\\
	&+\biggl\| \sum_{|\beta|=2} \int_0^{t/2} \int_{\mathbb{R}^3} \int_0^1 \frac{\nabla^\beta (\boldsymbol{a}\cdot\nabla) G(t, x-\lambda y)}{\beta!} \lambda y^\beta (u^2 - M_0^2G^2) (s,y) d\lambda dyds \biggr\|_{L^q (\mathbb{R}^3)}\\
	&\le
	\int_0^{t/2} \int_0^1 \| \partial_t (\boldsymbol{a}\cdot\nabla) G(t-\lambda s) \|_{L^q (\mathbb{R}^3)} s \| (u^2-M_0^2G^2) (s) \|_{L^1 (\mathbb{R}^3)} d\lambda ds\\
	&+ C \sum_{|\beta|=2} \int_0^{t/2} \| \nabla^\beta (\boldsymbol{a}\cdot\nabla) G(t) \|_{L^q (\mathbb{R}^3)} \| y^\beta (u^2-M_0^2G^2) (s) \|_{L^1 (\mathbb{R}^3)} ds\\
	&\le
	C \int_0^{t/2} \left( (t-s)^{-\gamma_q-\frac32} + t^{-\gamma_q-\frac32} \right) s^{-1/2} (1+s)^{-1/2} ds
	\le
	C t^{-\gamma_q-\frac32} \log (2+t)
\end{split}
\]
and
\[
\begin{split}
	&\biggl\| \int_{t/2}^t \int_{\mathbb{R}^3} \biggl( \boldsymbol{a}\cdot\nabla G(t-s,x-y) - \sum_{|\beta|=0}^1 \nabla^\beta (\boldsymbol{a}\cdot\nabla) G(t,x) (-y)^\beta \biggr)\\
	&\hspace{25mm} (u^2-M_0^2G^2) (s,y) dyds \biggr\|_{L^q (\mathbb{R}^3)}\\
	&\le
	\int_{t/2}^t \| \boldsymbol{a}\cdot\nabla G(t-s) \|_{L^1 (\mathbb{R}^3)} \| (u^2-M_0^2G^2) (s) \|_{L^q (\mathbb{R}^3)} ds\\
	&+
	\sum_{|\beta|=0}^1 \| \nabla^\beta (\boldsymbol{a}\cdot\nabla) G(t) \|_{L^q (\mathbb{R}^3)} \int_{t/2}^t \|  y^\beta (u^2-M_0^2G^2) (s) \|_{L^1 (\mathbb{R}^3)} ds
	\le
	Ct^{-\gamma_q-\frac32}
\end{split}
\]
from \eqref{EZexp} and \eqref{decay-wt}.
Hence, $\| r_3^1 (t) \|_{L^q (\mathbb{R}^3)} = O (t^{-\gamma_q-\frac32} \log t)$ as $t \to +\infty$.
We shape the first term of \eqref{bs13}.
Since
$
	\mathcal{F} [\boldsymbol{a}\cdot\nabla G(t-s)*G^2 (s)]
	=
	(4\pi)^{-3} s^{-3/2} \boldsymbol{a}\cdot i\xi e^{-(t-\frac{s}2)|\xi|^2},
$
we see that
\[
	\boldsymbol{a}\cdot\nabla G(t-s) * G^2 (s) = \frac{\sqrt{2\pi}}{32\pi^2} s^{-3/2} \boldsymbol{a}\cdot\nabla G(t-\tfrac{s}2).
\]
This formulation together with $\int_{\mathbb{R}^3} G^2 (s,y) dy = \frac{\sqrt{2\pi}}{32\pi^2} s^{-3/2}$ and $\int_{\mathbb{R}^3} y^\beta G^2 dy = 0$ for $|\beta| = 1$ yields that
\[
\begin{split}
	J_2 (t)
	&=
	\frac{M_0^2\sqrt{2\pi}}{32\pi^2} \int_0^t s^{-\frac32} \left( \boldsymbol{a}\cdot\nabla G(t-\tfrac{s}2) - \boldsymbol{a}\cdot\nabla G(t) \right) ds.
\end{split}
\]
Here, the mean value theorem guarantees integrability of this function as $s \to +0$.
Thus, by integration by parts, we have that
\[
\begin{split}
	&\int_0^t s^{-3/2} (\boldsymbol{a}\cdot\nabla G(t-\tfrac{s}2) - \boldsymbol{a}\cdot\nabla G(t)) ds\\
	&=
	-2 t^{-1/2} (\boldsymbol{a}\cdot\nabla G(\tfrac{t}2) - \boldsymbol{a}\cdot\nabla G(t))
	-
	\int_0^t s^{-1/2} \Delta (\boldsymbol{a}\cdot\nabla)G(t-\tfrac{s}2) ds.
\end{split}
\]
Substituting this result into \eqref{bs13}, we see for the second part of \eqref{bs1} that
\begin{equation}\label{bs13p}
\begin{split}
	&\int_0^t \int_{\mathbb{R}^3} \biggl( \boldsymbol{a}\cdot\nabla G(t-s,x-y) - \sum_{|\beta|=0}^1 \nabla^\beta (\boldsymbol{a}\cdot\nabla) G(t,x) (-y)^\beta \biggr) u^2 (s,y) dyds\\
	&=
	-\frac{\sqrt{2\pi}}{32\pi^2} M_0^2 \biggl( 2t^{-1/2} \left( \boldsymbol{a}\cdot\nabla G(\tfrac{t}2) - \boldsymbol{a}\cdot\nabla G(t) \right) + \int_0^t s^{-1/2} \Delta (\boldsymbol{a}\cdot\nabla) G(t-\tfrac{s}2) ds \biggr)\\
	&+ r_3^1 (t).
\end{split}
\end{equation}
The coefficient of fourth term of \eqref{bs1} is expanded as
\begin{equation}\label{bs1-3}
\begin{split}
	&- \sum_{|\beta|=0}^1 \nabla^\beta (\boldsymbol{a}\cdot\nabla) G(t) \int_t^\infty \int_{\mathbb{R}^3} (-y)^\beta u^2 (s,y) dyds\\
	&=
	- M_0^2 \boldsymbol{a}\cdot\nabla G(t) \int_t^\infty \int_{\mathbb{R}^3} G^2 (s,y) dyds + r_3^2 (t)\\
	&=
	- \frac{M_0^2 \sqrt{2\pi}}{16\pi^2} t^{-1/2} \boldsymbol{a}\cdot\nabla G(t) + r_3^2 (t)
\end{split}
\end{equation}
for
\begin{equation}\label{r32}
	r_3^2 (t) =- \sum_{|\beta|=0}^1 \nabla^\beta (\boldsymbol{a}\cdot\nabla) G(t) \int_t^\infty \int_{\mathbb{R}^3} (-y)^\beta (u^2-M_0^2G^2) (s,y) dyds
\end{equation}
since $\int_{\mathbb{R}^3} y^\beta G^2 dy = 0$ for $|\beta| = 1$.
Here, we putted
\[
	\int_t^\infty \int_{\mathbb{R}^3} G^2 (s,y) dyds
	= \frac{\sqrt{2\pi}}{16\pi^2} t^{-1/2}.
\]
This $r_3^2$ also plays similar role as one in Section \ref{subsec4dim}.
Indeed, the estimates  \eqref{EZexp} and \eqref{decay-wt} show
\[
	\biggl|  \int_t^\infty \int_{\mathbb{R}^3} (-y)^\beta (u^2-M_0^2G^2) (s,y) dyds \biggr|
	\le
	C \int_t^\infty s^{-2+\frac{|\beta|}2} ds
	=
	C t^{-1+\frac{|\beta|}2}
\]
and then $\| r_3^2 (t) \|_{L^q (\mathbb{R}^3)} = O (t^{-\gamma_q-\frac32})$ as $t \to +\infty$.
Adding \eqref{bs13p} and \eqref{bs1-3}, we see for the third and fourth parts of \eqref{bs1}  that
\[
\begin{split}
	&- \sum_{|\beta|=0}^1 \nabla^\beta (\boldsymbol{a}\cdot\nabla) G(t) \int_t^\infty \int_{\mathbb{R}^3} (-y)^\beta u^2 (s,y) dyds\\
	&+ \int_0^t \int_{\mathbb{R}^3} \biggl( \boldsymbol{a}\cdot\nabla G(t-s,x-y) - \sum_{|\beta|=0}^1 \nabla^\beta (\boldsymbol{a}\cdot\nabla) G(t,x) (-y)^\beta \biggr) u^2 (s,y) dyds\\
	&=  - \frac{\sqrt{2\pi}}{32\pi^2} M_0^2 \left( 2t^{-1/2} \boldsymbol{a}\cdot\nabla G(\tfrac{t}2) + \int_0^t s^{-1/2} \Delta (\boldsymbol{a}\cdot\nabla) G(t-\tfrac{s}2) ds \right) + r_3^1 (t) + r_3^2 (t).
\end{split}
\]
Now we see the remaind part of $U_2$.
Therefore, by substituting this to \eqref{bs1}, we obtain that $u(t) = U_0 (t) + U_1 (t) + U_2 (t) + r_3 (t)$ for $U_m$ given by \eqref{U0U1} and \eqref{U2}, and $r_3 = r_3^0 + r_3^1 + r_3^2$.
We already confirmed that $\| r_3 (t) \|_{L^q (\mathbb{R}^3)} = O(t^{-\gamma_q-\frac32}\log t)$ and then
\begin{equation}\label{EZexpK}
	\| u(t) - M_0 G (t) - \boldsymbol{M}_1 \cdot \nabla G (t) - U_2 (t) \|_{L^q (\mathbb{R}^3)} = O(t^{-\gamma_q-\frac32}\log t)
\end{equation}
as $t \to +\infty$.
The logarithmic evolution in \eqref{EZexp} is eliminated.
\subsection{The expansion up to third order}
In the last section, we saw that $u = U_0 + U_1 + U_2 + r_3$ and $U_m$ has no logarithms.
We expand the first term of $r_3$ defined by \eqref{r31} as
\[
\begin{split}
	&r_3^1 (t)
	=
	\sum_{2l+|\beta|=2} \frac{\partial_t \nabla^\beta (\boldsymbol{a}\cdot\nabla) G(t)}{\beta!} \int_0^t \int_{\mathbb{R}^3}
		(-s)^l (-y)^\beta (u^2 - M_0^2 G^2) (s,y)
	dyds\\
	&+
	\int_0^t \int_{\mathbb{R}^3} \biggl(
		\boldsymbol{a}\cdot\nabla G(t-s,x-y) - \sum_{2l+|\beta|=0}^2 \frac{\partial_t^l \nabla^\beta (\boldsymbol{a}\cdot\nabla) G(t,x)}{\beta!} (-s)^l (-y)^\beta \biggr)\\
		&\hspace{20mm} (u^2 - M_0^2 G^2) (s,y)
	dyds.
\end{split}
\]
In the renormalization process, we remark that $\int_{\mathbb{R}^3} (-s)^l (-y)^\beta (\boldsymbol{M}_1\cdot\nabla) (G^2) dy = 0$ for $2l+|\beta|=2$.
Therefore, the renormalization expands these terms as
\[
\begin{split}
	&\int_0^t \int_{\mathbb{R}^3}
		(-s)^l (-y)^\beta (u^2 - M_0^2 G^2) (s,y)
	dyds\\
	&=
	\int_0^\infty \int_{\mathbb{R}^3}
		(-s)^l (-y)^\beta (u^2 - M_0^2 G^2) (s,y)
	dyds\\
	&- \int_t^\infty \int_{\mathbb{R}^3}
		(-s)^l (-y)^\beta (u^2 - M_0^2 G^2 - M_0 \boldsymbol{M}_1 \cdot \nabla (G^2)) (s,y)
	dyds
\end{split}
\]
and then
\begin{equation}\label{r31bs}
\begin{split}
	&r_3^1 (t)
	=
	\sum_{2l+|\beta|=2} \frac{\partial_t \nabla^\beta (\boldsymbol{a}\cdot\nabla) G(t)}{\beta!} \int_0^\infty \int_{\mathbb{R}^3}
		(-s)^l (-y)^\beta (u^2 - M_0^2 G^2) (s,y)
	dyds + J_3 (t)\\
	&- \sum_{2l+|\beta|=2} \frac{\partial_t \nabla^\beta (\boldsymbol{a}\cdot\nabla) G(t)}{\beta!} \int_t^\infty \int_{\mathbb{R}^3}
		(-s)^l (-y)^\beta (u^2 - M_0^2 G^2\\
		&\hspace{30mm} - M_0 \boldsymbol{M}_1 \cdot \nabla (G^2)) (s,y)
	dyds\\
	& + r_4^1 (t)
\end{split}
\end{equation}
for
\[
\begin{split}
	&J_3 (t) = M_0\int_0^t \int_{\mathbb{R}^3} \biggl(
		\boldsymbol{a}\cdot\nabla G(t-s,x-y) - \sum_{2l+|\beta|=0}^2 \frac{\partial_t^l \nabla^\beta (\boldsymbol{a}\cdot\nabla) G(t,x)}{\beta!} (-s)^l (-y)^\beta \biggr)\\
		&\hspace{30mm} (\boldsymbol{M}_1\cdot\nabla) (G^2) (s,y)
	dyds
\end{split}
\]
and
\begin{equation}\label{r41}
\begin{split}
	&r_4^1 (t)
	=
	\int_0^t \int_{\mathbb{R}^3} \biggl(
		\boldsymbol{a}\cdot\nabla G(t-s,x-y) - \sum_{2l+|\beta|=0}^2 \frac{\partial_t^l \nabla^\beta (\boldsymbol{a}\cdot\nabla) G(t,x)}{\beta!} (-s)^l (-y)^\beta \biggr)\\
		&\hspace{40mm} (u^2 - M_0^2 G^2 - M_0 \boldsymbol{M}_1\cdot\nabla (G^2)) (s,y)
	dyds.
\end{split}
\end{equation}
Here, integrability of the coefficients on the first part of \eqref{r31bs} is guaranteed by \eqref{EZexp} and \eqref{decay-wt} since $\int_{\mathbb{R}^3} y^\beta (u^2 - M_0^2 G^2) dy = \int_{\mathbb{R}^3} y^\beta (u^2 - M_0^2 G^2 - M_0 \boldsymbol{M}_1 \cdot \nabla (G^2)) dy$ for $2l+|\beta| = 2$ because $y^\beta \boldsymbol{M}_1 \cdot \nabla (G^2)$ should be odd in some variable.
The similar estimate as in Section \ref{subsect3dim1} together with \eqref{EZexpK} yields that
$
	\| r_4^1 (t) \|_{L^q (\mathbb{R}^3)}
	\le Ct^{-\gamma_q-2} \log (2+t).
$
We shape $J_3$.
Since
\[
	\mathcal{F} [\boldsymbol{a}\cdot \nabla G(t-s)* (\boldsymbol{M}_1\cdot\nabla G^2) (s)]
	=
	\frac{\sqrt{2\pi}}{32\pi^2} s^{-3/2} (2\pi)^{-3/2} (\boldsymbol{M}_1 \cdot i\xi) (\boldsymbol{a}\cdot i\xi) e^{-(t-\frac{s}2)|\xi|^2},
\]
we have for the former part that
\[
	\boldsymbol{a}\cdot \nabla G(t-s)* (\boldsymbol{M}_1\cdot\nabla G^2) (s)
	= \frac{\sqrt{2\pi}}{32\pi^2} s^{-3/2} (\boldsymbol{M}_1\cdot\nabla)(\boldsymbol{a}\cdot\nabla) G(t-\tfrac{s}2).
\]
On the other hand, if we omit the integrals that clearly disappear, the latter part is rewritten as
\[
\begin{split}
	&-\sum_{2l+|\beta|=0}^2 \frac{\partial_t^l \nabla^\beta (\boldsymbol{a}\cdot\nabla) G(t)}{\beta!} (-s)^l \int_{\mathbb{R}^3} (-y)^\beta (\boldsymbol{M}_1\cdot\nabla) (G^2) (s,y) dy\\
	&=
	\sum_{j=1}^3 \partial_j (\boldsymbol{a}\cdot\nabla) G(t) \int_{\mathbb{R}^3} y_j M_1^j \partial_j (G^2) (s,y) dy
	=
	- \frac{\sqrt{2\pi}}{32\pi^2} s^{-3/2} (\boldsymbol{M}_1\cdot\nabla)(\boldsymbol{a}\cdot\nabla) G(t).
\end{split}
\]
Hence, $J_3$ on \eqref{r31bs} is rewritten as
\[
\begin{split}
	J_3 (t)
	&=
	\frac{M_0 \sqrt{2\pi}}{32\pi^2} \int_0^t s^{-3/2} \left( (\boldsymbol{M}_1\cdot\nabla) (\boldsymbol{a}\cdot\nabla) G(t-\tfrac{s}2) - (\boldsymbol{M}_1\cdot\nabla) (\boldsymbol{a}\cdot\nabla) G(t) \right) ds.
\end{split}
\]
Here, the mean value theorem mitigates the singularity as $s \to + 0$.
Thus, the integration by parts provides that
\[
\begin{split}
	&\int_0^t s^{-3/2} \left( (\boldsymbol{M}_1\cdot\nabla) (\boldsymbol{a}\cdot\nabla) G(t-\tfrac{s}2) - (\boldsymbol{M}_1\cdot\nabla) (\boldsymbol{a}\cdot\nabla) G(t) \right) ds\\
	&=
	- 2t^{-1/2} \left( (\boldsymbol{M}_1\cdot\nabla)(\boldsymbol{a}\cdot\nabla) G(\tfrac{t}2) - (\boldsymbol{M}_1\cdot\nabla)(\boldsymbol{a}\cdot\nabla) G (t) \right)\\
	&- \int_0^t s^{-1/2} \Delta (\boldsymbol{M}_1\cdot\nabla) (\boldsymbol{a}\cdot\nabla) G(t-\tfrac{s}2) ds.
\end{split}
\]
To summarize them, we see
\begin{equation}\label{r31fin}
\begin{split}
	&r_3^1 (t)
	=
	\sum_{2l+|\beta|=2} \frac{\partial_t \nabla^\beta (\boldsymbol{a}\cdot\nabla) G(t)}{\beta!} \int_0^\infty \int_{\mathbb{R}^3}
		(-s)^l (-y)^\beta (u^2 - M_0^2 G^2) (s,y)
	dyds\\
	&
	-\frac{M_0\sqrt{2\pi}}{16\pi^2} t^{-1/2} \left( (\boldsymbol{M}_1\cdot\nabla)(\boldsymbol{a}\cdot\nabla) G(\tfrac{t}2) - (\boldsymbol{M}_1\cdot\nabla) (\boldsymbol{a}\cdot\nabla) G(t) \right)\\
	&- \frac{M_0\sqrt{2\pi}}{32\pi^2}\int_0^t s^{-1/2} \Delta (\boldsymbol{M}_1\cdot\nabla)(\boldsymbol{a}\cdot\nabla) G(t-\tfrac{s}2) ds\\
	&- \sum_{2l+|\beta|=2} \frac{\partial_t \nabla^\beta (\boldsymbol{a}\cdot\nabla) G(t)}{\beta!} \int_t^\infty \int_{\mathbb{R}^3}
		(-s)^l (-y)^\beta (u^2 - M_0^2 G^2\\
		&\hspace{30mm} - M_0 \boldsymbol{M}_1 \cdot \nabla (G^2)) (s,y)
	dyds\\
	&+ r_4^1 (t).
\end{split}
\end{equation}
Next, we reform $r_3^2$ given by \eqref{r32}.
We renormalize $u^2 - M_0^2 G^2$ in $r_3^2$ by $2U_0 U_1 = M_0 \boldsymbol{M}_1 \cdot \nabla (G^2)$, then, since $\int_{\mathbb{R}^3} \nabla (G^2) dy = 0$, we see that
\[
\begin{split}
	&r_3^2 (t)
	= M_0 \sum_{|\beta| = 1} \nabla^\beta (\boldsymbol{a}\cdot\nabla) G(t) \int_t^\infty \int_{\mathbb{R}^3} y^\beta \boldsymbol{M}_1\cdot\nabla (G^2) (s,y) dyds\\
	&- \sum_{|\beta|=0}^1 \nabla^\beta (\boldsymbol{a}\cdot\nabla) G(t) \int_t^\infty \int_{\mathbb{R}^3}  (-y)^\beta (u^2 - M_0^2 G^2 - M_0 \boldsymbol{M}_1 \cdot \nabla (G^2)) (s,y)  dyds.
\end{split}
\]
If we denote $\boldsymbol{M}_1 = (M_1^1,M_1^2,M_1^3)$, then we can rewrite the first term by
\[
\begin{split}
	&M_0\sum_{|\beta| = 1} \nabla^\beta (\boldsymbol{a}\cdot\nabla) G(t) \int_t^\infty \int_{\mathbb{R}^3} y^\beta \boldsymbol{M}_1\cdot\nabla (G^2) (s,y) dyds\\
	&=
	M_0 \sum_{j=1}^3 \partial_j (\boldsymbol{a}\cdot\nabla) G(t) \int_t^\infty s^{-3/2} ds \int_{\mathbb{R}^3} y_j M_1^j \partial_j (G^2) (1,y) dyds\\
	&= -\frac{\sqrt{2\pi}}{16\pi^2}M_0 t^{-1/2} (\boldsymbol{M}_1\cdot\nabla)(\boldsymbol{a}\cdot\nabla) G(t)
\end{split}
\]
since $\int_{\mathbb{R}^3} y_k \partial_j (G^2) dy = 0$ for $j \neq k$.
Substituting this result into $r_3^2$ and adding it to \eqref{r31fin} provide
\[
\begin{split}
	&r_3^1 (t) + r_3^2 (t)
	=
	\sum_{2l+|\beta|=2} \frac{\partial_t^l \nabla^\beta (\boldsymbol{a}\cdot\nabla) G(t)}{\beta!} \int_0^\infty \int_{\mathbb{R}^3}
		(-s)^l (-y)^\beta (u^2 - M_0^2G^2) (s,y)
	dyds\\
	&- \frac{M_0\sqrt{2\pi}}{16\pi^2} t^{-1/2} (\boldsymbol{M}_1\cdot\nabla)(\boldsymbol{a}\cdot\nabla) G(\tfrac{t}2)- \frac{M_0\sqrt{2\pi}}{32\pi^2} \int_0^t s^{-1/2} \Delta (\boldsymbol{M}_1\cdot\nabla)(\boldsymbol{a}\cdot\nabla) G(t-\tfrac{s}2) ds\\
	&+r_4^1 (t) + r_4^2 (t)
\end{split}
\]
for
\[
\begin{split}
	&r_4^2 (t)= - \sum_{2l+|\beta|=0}^2 \frac{\partial_t^l \nabla^\beta (\boldsymbol{a}\cdot\nabla) G(t)}{\beta!} \int_t^\infty \int_{\mathbb{R}^3} (-s)^l (-y)^\beta (u^2 - M_0^2 G^2\\
	&\hspace{30mm}  - M_0 \boldsymbol{M}_1 \cdot \nabla (G^2)) (s,y)  dyds.
\end{split}
\]
Now we find the part of $U_3$ introduced by \eqref{U3}.
By employing \eqref{decay-wt} and \eqref{EZexpK}, we see for $r_4^2$ that
\[
\begin{split}
	&\biggl| \int_t^\infty \int_{\mathbb{R}^3}  (-y)^\beta (u^2 - M_0^2 G^2 - M_0 \boldsymbol{M}_1 \cdot \nabla (G^2)) (s,y)  dyds \biggr|
	\le
	C t^{-1+\frac{|\beta|}2}.
\end{split}
\]
Hence, $\| r_4^2 (t) \|_{L^q (\mathbb{R}^3)} \le Ct^{-\gamma_q-2}$.
Therefore, we obtain that $r_3 = U_3 + r_4$ and then $u = U_0 + U_1 + U_2 + U_3 + r_4$ for
$r_4 = r_4^0 + r_4^1 + r_4^2$.
Of coursely, $r_4^0$ is given by \eqref{rm0}.
We already checked that $\| r_4 (t) \|_{L^q (\mathbb{R}^3)} = O (t^{-\gamma_q-2}\log t)$ as $t \to +\infty$.
Consequently, we complete the asymptotic expansion up to third order.
\subsection{Derivation of the logarithmic evolution}
So far, we expand the solution as $u = U_0 + U_1 + U_2 + U_3 + r_4$ and $U_m$ contains no logarithms.
Here, $U_m$ are given by \eqref{U0U1}-\eqref{U3}.
The logarithmic evolution may be hyden in $r_4 = r_4^0 + r_4^1 + r_4^2$.
As we already checked, $r_4^0$ and $r_4^2$ contain no logarithms.
We expand $r_4^1$ given by \eqref{r41} and renormalize $u^2 - M_0^2 G^2 - M_0 (\boldsymbol{M}_1\cdot\nabla) G^2 = u^2 - U_0^2 - 2U_0 U_1$ by $(\boldsymbol{M}_1\cdot \nabla G)^2+ 2M_0 G U_2 = U_1^2 + 2U_0 U_2$, then we obtain
\begin{equation}\label{bs14}
\begin{split}
	&r_4^1 (t)\\
	& =
	\sum_{2l+|\beta|=3} \frac{\partial_t^l \nabla^\beta (\boldsymbol{a}\cdot\nabla) G(t)}{\beta!} \int_0^t \int_{\mathbb{R}^3}
		(-s)^l (-y)^\beta \left( (\boldsymbol{M}_1 \cdot G)^2 + 2M_0 GU_2 \right) (1+s,y)
	dyds\\
	&+ r_4^3 (t) + r_4^4 (t)
\end{split}
\end{equation}
for
\[
\begin{split}
	&r_4^3 (t)\\
	& = \sum_{2l+|\beta|=3} \frac{\partial_t^l \nabla^\beta (\boldsymbol{a}\cdot\nabla) G(t)}{\beta!} \int_0^t \int_{\mathbb{R}^3}
		(-s)^l (-y)^\beta\bigl( (u^2 - M_0^2 G^2 - M_0\boldsymbol{M}_1\cdot\nabla (G^2)) (s,y)\\
		&\hspace{15mm}- ((\boldsymbol{M}_1\cdot\nabla G)^2 + M_0 GU_2) (1+s,y) \bigr)
	dyds
\end{split}
\]
and
\[
\begin{split}
	&r_4^4 (t)
	= \int_0^t \int_{\mathbb{R}^3} \biggl(
		\boldsymbol{a}\cdot\nabla G(t-s,x-y) - \sum_{2l+|\beta|=0}^3 \frac{\partial_t^l \nabla^\beta (\boldsymbol{a}\cdot\nabla) G(t,x)}{\beta!} (-s)^l (-y)^\beta \biggr)\\
		&\hspace{40mm} (u^2 - M_0^2 G^2 - M_0\boldsymbol{M}_1\cdot\nabla (G^2)) (s,y)
	dyds.
\end{split}
\]
A coupling of \eqref{decay-wt} and \eqref{EZexpK} yields that
\[
\begin{split}
	&\biggl| \int_0^t \int_{\mathbb{R}^3}
		(-s)^l (-y)^\beta\bigl( (u^2 - M_0^2 G^2 -  M_0\boldsymbol{M}_1\cdot\nabla (G^2) ) (s,y)\\
		&\hspace{15mm}- ((\boldsymbol{M}_1\cdot\nabla G)^2 + M_0 GU_2) (1+s,y) \bigr) dyds
	dyds \biggr|\\
	&\le
	C \int_0^t (1+s)^{-\frac32} \log (2+s) ds < + \infty,
\end{split}
\]
then $\| r_4^3 (t) \|_{L^q (\mathbb{R}^3)} \le Ct^{-\gamma_q-2}$.
This part contains no logarithms.
The similar argument as in Section \ref{subsec4dim} together with \eqref{decay-wt} and \eqref{EZexpK} gives that
\[
\begin{split}
	&\| r_4^4 (t) \|_{L^q (\mathbb{R}^3)}
	\le
	\frac12 \int_0^{t/2} \int_0^1
		\| \partial_t^2 (\boldsymbol{a}\cdot\nabla) G(t-\lambda s) \|_{L^q (\mathbb{R}^3)}\\
		&\hspace{15mm}
		s^2 \|  (u^2 - M_0^2 G^2 -  M_0\boldsymbol{M}_1\cdot\nabla (G^2)) (s) \|_{L^1 (\mathbb{R}^3)}
	d\lambda ds\\
	&+C \sum_{2l+|\beta|=4} \int_0^{t/2}
		\| \partial_t^l \nabla^\beta (\boldsymbol{a}\cdot\nabla) G(t) \|_{L^q (\mathbb{R}^3)}\\
		&\hspace{15mm}
		s^l \| y^\beta (u^2 - M_0^2 G^2 -  M_0\boldsymbol{M}_1\cdot\nabla (G^2)) (s) \|_{L^1 (\mathbb{R}^3)}
	ds\\
	&+ \int_{t/2}^t \| \boldsymbol{a}\cdot\nabla G(t-s) \|_{L^1 (\mathbb{R}^3)} \| (u^2 - M_0^2 G^2 -  M_0\boldsymbol{M}_1\cdot\nabla (G^2)) (s) \|_{L^q (\mathbb{R}^3)}ds\\
	&+ C\sum_{2l+|\beta|=0}^3 \| \partial_t^l \nabla^\beta (\boldsymbol{a}\cdot\nabla) G(t) \|_{L^q (\mathbb{R}^3)} \int_{t/2}^t s^l \| y^\beta (u^2 - M_0^2 G^2\\
	&\hspace{50mm} -  M_0\boldsymbol{M}_1\cdot\nabla (G^2)) (s) \|_{L^1 (\mathbb{R}^3)}
	ds\\
	&= O(t^{-\gamma_q-2})
\end{split}
\]
as $t \to +\infty$.
This part does not provide logarithms too.
The logarithmic evolution could be included in the former part of \eqref{bs14} since
\[
\begin{split}
	&\int_0^t \int_{\mathbb{R}^3}
		(-s)^l (-y)^\beta \left( (\boldsymbol{M}_1 \cdot G)^2 + 2M_0 GU_2 \right) (1+s,y)
	dyds\\
	&=
	\int_0^t s^l (1+s)^{-l-1} ds \int_{\mathbb{R}^3} (-1)^l (-y)^\beta \left( (\boldsymbol{M}_1 \cdot \nabla G)^2 + 2M_0 GU_2 \right) (1,y)
	dy.
\end{split}
\]
Noting that $2l+|\beta|=3$, several integrands are vanishing.
Clearly, $\int_{\mathbb{R}^3} (-1)^l (-y)^\beta (\boldsymbol{M}_1\cdot\nabla G)^2 (1,y) dy = 0$.
We separate $U_2$ by $U_2 = U_2^{\mathrm{evn}} + U_2^{\mathrm{odd}}$ for
\[
\begin{split}
	&U_2^{\mathrm{evn}} (t) = \sum_{|\alpha|=2} \frac{\nabla^\alpha G(t)}{\alpha!} \int_{\mathbb{R}^3} y^\alpha u_0 (y) dy
	- \nabla (\boldsymbol{a}\cdot\nabla) G(t) \cdot \int_0^\infty \int_{\mathbb{R}^3} y u^2 (s,y) dyds
\end{split}
\]
and
\[
\begin{split}
	&U_2^{\mathrm{odd}} (t) =-  \frac{\sqrt{2\pi}}{32\pi^2} M_0^2 \left( 2t^{-1/2} \boldsymbol{a}\cdot\nabla G(\tfrac{t}2) + \int_0^t s^{-1/2} \Delta (\boldsymbol{a}\cdot\nabla) G(t-\tfrac{s}2) ds \right),
\end{split}
\]
then $\int_{\mathbb{R}^3} (-1)^l (-y)^\beta (2M_0 GU_2^{\mathrm{evn}}) (1,y) dy = 0$ since $2l+|\beta| = 3$ and $(GU_2^{\mathrm{evn}}) (1,-y) = (GU_2^{\mathrm{evn}}) (1,y)$.
Thus, the coefficient of former part of \eqref{bs14} is simplified as
\[
\begin{split}
	&\int_0^t \int_{\mathbb{R}^3}
		(-s)^l (-y)^\beta \left( (\boldsymbol{M}_1 \cdot G)^2 + 2M_0 GU_2 \right) (1+s,y)
	dyds\\
	&=
	2 M_0 \int_0^t s^l (1+s)^{-l-1} ds \int_{\mathbb{R}^3} (-1)^l (-y)^\beta (GU_2^{\mathrm{odd}}) (1,y) dy,
\end{split}
\]
and
\[
\begin{split}
	&\int_{\mathbb{R}^3} (-1)^l (-y)^\beta (GU_2^{\mathrm{odd}}) (1,y) dy\\
	&=
	- \frac{\sqrt{2\pi}}{32\pi^2} M_0^2 \int_{\mathbb{R}^3} (-1)^l (-y)^\beta G(1,y) \biggl( 2\boldsymbol{a}\cdot\nabla G(\tfrac12,y)\\
	&\hspace{20mm}  + \int_0^1 s^{-1/2} \Delta (\boldsymbol{a}\cdot\nabla) G(1-\tfrac{s}2,y) ds \biggr) dy.
\end{split}
\]
Here, $\int_0^t s^l (1+s)^{-l-1} ds$ leads the desired logarithmic evolution, that is, $\int_0^t s^l (1+s)^{-l-1} ds = \log t + O (1)$ as $t \to +\infty$.
Summarizing what we have discussed so far, $r_4 = K_4 \log t + \tilde{r}_4$ and then the solution $u$ is expanded as
\[
	u(t) = U_0 (t) + U_1 (t) + U_2 (t) + U_3 (t) + K_4 (t) \log t + \tilde{r}_4 (t)
\]
for $U_m$ given by \eqref{U0U1}-\eqref{U3}, and
\[
\begin{split}
	&K_4 (t) = - \frac{\sqrt{2\pi}}{16\pi^2} M_0^3 \sum_{2l+|\beta|=3} \frac{\partial_t^l \nabla^\beta (\boldsymbol{a}\cdot\nabla) G(t)}{\beta!} \int_{\mathbb{R}^3} (-1)^l (-y)^\beta G(1,y)\\
	&\hspace{15mm} \left( 2\boldsymbol{a}\cdot\nabla G(\tfrac12,y) + \int_0^1 s^{-1/2} \Delta (\boldsymbol{a}\cdot\nabla) G(1-\tfrac{s}2,y) ds \right) dy,
\end{split}
\]
and some $\tilde{r}_4$ satisfying $\| \tilde{r}_4 (t) \|_{L^q (\mathbb{R}^3)} = O (t^{-\gamma_q-2})$ as $t\to +\infty$.
Here, $\tilde{r}_4$ has a messy form, but it can be written down as $\tilde{r}_4 = r_4^0 + r_4^2 + r_4^3 + r_4^4 + r_4^5$ for the above errors and
\[
\begin{split}
	&r_4^5 (t) = - \frac{\sqrt{2\pi}}{16\pi^2} M_0^3 \sum_{2l+|\beta|=3} \frac{\partial_t^l \nabla^\beta (\boldsymbol{a}\cdot\nabla) G(t)}{\beta!} \left( \int_0^t s^l (1+s)^{-l-1} ds - \log t \right)\\
	&\hspace{10mm}\int_{\mathbb{R}^3} (-1)^l (-y)^\beta G(1,y)
	\left( 2\boldsymbol{a}\cdot\nabla G(\tfrac12,y) + \int_0^1 s^{-1/2} \Delta (\boldsymbol{a}\cdot\nabla) G(1-\tfrac{s}2,y) ds \right) dy
\end{split}
\]
which has the same root as $K_4 \log t$.
We should show that this $K_4$ is written as \eqref{K4}.
At this stage, $K_4$ contains extra integrands yet.
Generally speaking that an integral of odd-type $f$, that is $f(-y) = - f(y)$, is vanishing.
We omit them and obtain that
\[
\begin{split}
	&K_4 (t) = - \frac{\sqrt{2\pi}}{16\pi^2} M_0^3 \sum_{j=1}^3 a_j \partial_t \partial_j (\boldsymbol{a}\cdot\nabla) G(t) \int_{\mathbb{R}^3} y_j G(1,y)
	\biggl( 2 \partial_j G(\tfrac12,y)\\
	&\hspace{50mm}  + \int_0^1 s^{-1/2} \Delta \partial_j G(1-\tfrac{s}2,y) ds \biggr) dy\\
	&+ \frac{\sqrt{2\pi}}{16\pi^2} M_0^3 \sum_{j=1}^3 \sum_{k\neq j} \frac{a_j\partial_k^2 \partial_j (\boldsymbol{a}\cdot\nabla) G(t)}{2!} \int_{\mathbb{R}^3} y_k^2 y_j G(1,y)
	\biggl( 2 \partial_j G(\tfrac12,y)\\
	&\hspace{50mm}  + \int_0^1 s^{-1/2} \Delta \partial_j G(1-\tfrac{s}2,y) ds \biggr) dy\\
	&+\frac{\sqrt{2\pi}}{16\pi^2} M_0^3 \sum_{j=1}^3 \frac{a_j \partial_j^3 (\boldsymbol{a}\cdot\nabla) G(t)}{3!} \int_{\mathbb{R}^3} y_j^3 G(1,y)
	\biggl( 2 \partial_j G(\tfrac12,y)\\
	&\hspace{50mm}  + \int_0^1 s^{-1/2} \Delta \partial_j G(1-\tfrac{s}2,y) ds \biggr) dy.
\end{split}
\]
We calculate these integrals.
The elementary calculus provide that $\int_{\mathbb{R}^3} y_j G(1,y) \partial_j G(\frac12,y) dy = - \frac{\sqrt{6\pi}}{54\pi^2},$ $\int_{\mathbb{R}^3} y_k^2 y_j G(1,y) \partial_j G(\frac12,y) dy = - \frac{\sqrt{6\pi}}{81\pi^2}$ for $k\neq j$, and $\int_{\mathbb{R}^3} y_j^3 G(1,y) \partial_j G(\frac12,y) dy = - \frac{\sqrt{6\pi}}{27\pi^2}$.
Plancherel theorem yields that
\[
\begin{split}
	&\int_{\mathbb{R}^3} y_j G(1,y) \int_0^1 s^{-1/2} \Delta \partial_j G(1-\tfrac{s}2,y) ds dy
	=
	\frac1{4\pi^3} \int_0^1 s^{-1/2} \int_{\mathbb{R}^3} \xi_j^2 |\xi|^2 e^{-(2-\frac{s}2)|\xi|^2} d\xi ds\\
	&=
	\frac1{4\pi^3} \int_0^1 s^{-1/2} (2-\tfrac{s}2)^{-7/2} ds \int_{\mathbb{R}^3} \xi_j^2 |\xi|^2 e^{-|\xi|^2} d\xi
	= \frac{7\sqrt{6\pi}}{2^3 \cdot 3^3 \pi^2}.
\end{split}
\]
Here we substituted as $s = 4\sin^2 \theta$ and then we saw
\[
	\int_0^1 s^{-1/2} (2-\tfrac{s}2)^{-7/2} ds = \frac{\sqrt{2}}4 \int_0^{\pi/6} \frac{d\theta}{\cos^6\theta} = \frac{14\sqrt{6}}{3^3 \cdot 5}.
\]
Similarly
\[
\begin{split}
	&\int_{\mathbb{R}^3} y_k^2 y_j G(1,y) \int_0^1 s^{-1/2} \Delta \partial_j G(1-\tfrac{s}2,y) ds dy
	=
	\frac{11\sqrt{6\pi}}{2 \cdot 3^4 \cdot 5\pi^2}
\end{split}
\]
for $k \neq j$, and
\[
\begin{split}
	&\int_{\mathbb{R}^3} y_j^3 G(1,y) \int_0^1 s^{-1/2} \Delta \partial_j G(1-\tfrac{s}2,y) ds dy
	=
	\frac{11\sqrt{6\pi}}{2\cdot 3^3 \cdot 5\pi^2}.
\end{split}
\]
Clearly, they are independent of $j$ or $k$.
Substituting these results into the last $K_4$ shows us \eqref{K4} and we complete the proof.

\section{The other odd-dimensional cases}\label{sect-oddim}
We show Proposition \ref{propoddim}.
The procedure is an extension of that described in Section \ref{sect3dim} and is similar as in \cite{Ym25}.
Specifically, we renormalize $u^2$ in the nonlinear term by the production of $U_m$ given by \eqref{Um}.
In this procedure, the spatial structures of these $U_m$ are important.
When $m$ is odd, $U_m$ is odd-type that is $U_m (t,-x) = - U_m (t,x)$ for $(t,x) \in \mathbb{R}_+ \times \mathbb{R}^n$.
Oppositely, this is even-type that is $U_m (t,-x) = U_m (t,x)$.
Consequently,  if $|\beta| + m_1 + m_2$ is odd, then
\begin{equation}\label{odd-van}
	\int_{\mathbb{R}^n} x^\beta (U_{m_1} U_{m_2}) (t,x) dx = 0
\end{equation}
for $t > 0$.
Since only these parities and the parabolic scales are important, we will not write the specific form of $U_m$ in this section.
Firstly, we expand the solution as
%
\begin{equation}\label{bsn1}
\begin{split}
	&u(t) = \sum_{|\alpha|=0}^{n-1} \int_{\mathbb{R}^n} (-y)^\alpha u_0 (y) dy\\
	&+ \sum_{2l+|\beta|=0}^{n-3} \frac{\partial_t^l \nabla^\beta (\boldsymbol{a}\cdot\nabla) G(t)}{l!\beta!} \int_0^t \int_{\mathbb{R}^n} (-s)^l (-y)^\beta u^2 (s,y) dyds\\
	&+ \sum_{2l+|\beta|=n-2} \frac{\partial_t^l \nabla^\beta (\boldsymbol{a}\cdot\nabla) G(t)}{l!\beta!} \int_0^t \int_{\mathbb{R}^n} (-s)^l (-y)^\beta (u^2-U_0^2) (s,y) dyds\\
	&+ \int_0^t \int_{\mathbb{R}^n} \biggl( \boldsymbol{a}\cdot \nabla G(t-s,x-y) - \sum_{2l+|\beta|=0}^{n-2} \frac{\partial_t^l \nabla^\beta (\boldsymbol{a}\cdot\nabla) G(t,x)}{l!\beta!} (-s)^l (-y)^\beta \biggr)\\
	&\hspace{20mm} u^2 (s,y) dyds
	+ r_n^0 (t)
\end{split}
\end{equation}
for $r_n^0$ given by \eqref{rm0}, where \eqref{odd-van} is applied in the third part.
%
Here, from \eqref{decay}, \eqref{EZexp} and \eqref{decay-wt}, the coefficients of the second part fulfill
\[
	\biggl| \int_{\mathbb{R}^n} (-s)^l (-y)^\beta u^2 (s,y) dy \biggr|
	\le
	C (1+s)^{-\frac{n}2+l+\frac{|\beta|}2}
\]
for $0 \le 2l+|\beta| \le n-3$, and one of the third part satisfies
\[
	\bigg| \int_{\mathbb{R}^n} (-s)^l (-y)^\beta (u^2-U_0^2) (s,y) dy \biggr|
	\le C (1+s)^{-3/2}
\]
for $2l+|\beta| = n-2$ since $U_0^2$ is a false term.
Hence, the coefficients are divided to
\[
\begin{split}
	&\int_0^t \int_{\mathbb{R}^n} (-s)^l (-y)^\beta u^2 (s,y) dyds\\
	&=
	\int_0^\infty \int_{\mathbb{R}^n} (-s)^l (-y)^\beta u^2 (s,y) dyds
	-
	\int_t^\infty \int_{\mathbb{R}^n} (-s)^l (-y)^\beta u^2 (s,y) dyds
\end{split}
\]
and
\[
	\biggl| \int_t^\infty \int_{\mathbb{R}^n} (-s)^l (-y)^\beta u^2 (s,y) dyds \biggr|
	\le
	C t^{-\frac{n}2+1+l+\frac{|\beta|}2}
\]
for $0 \le 2l+|\beta| \le n-3$, and
\[
\begin{split}
	&\int_0^t \int_{\mathbb{R}^n} (-s)^l (-y)^\beta (u^2-U_0^2) (s,y) dyds\\
	&=
	\int_0^\infty \int_{\mathbb{R}^n} (-s)^l (-y)^\beta(u^2-U_0^2) (s,y) dyds
	-
	\int_t^\infty \int_{\mathbb{R}^n} (-s)^l (-y)^\beta (u^2-U_0^2) (s,y) dyds
\end{split}
\]
and
\[
	\biggl| \int_t^\infty \int_{\mathbb{R}^n} (-s)^l (-y)^\beta (u^2-U_0^2) (s,y) dyds \biggr|
	\le
	C t^{-1/2}
\]
for $2l+|\beta| = n-2$.
Here the treatment for $2l+|\beta| = n-2$ is a novelty.
Taylor theorem together with \eqref{decay} and \eqref{decay-wt} guarantees for the fourth term of \eqref{bsn1} that
\[
\begin{split}
	&\biggl\| \int_0^t \int_{\mathbb{R}^n} \biggl( \boldsymbol{a}\cdot\nabla G(t-s,x-y) - \sum_{2l+|\beta|=0}^{n-2} \frac{\partial_t^l \nabla^\beta (\boldsymbol{a}\cdot\nabla) G(t,x)}{l!\beta!} (-s)^l (-y)^\beta \biggr)\\
	&\hspace{20mm} u^2 (s,y) dyds \biggr\|_{L^q (\mathbb{R}^n)}
	= O (t^{-\gamma_q-\frac{n}2 + \frac12})
\end{split}
\]
as $t \to + \infty$.
At this stage, we confirm that the logarithmic evolution in \eqref{EZexp} is false and then
\begin{equation}\label{EZexpK2}
	\biggl\| u(t) - \sum_{m=0}^{n-2} U_m (t) \biggr\|_{L^q (\mathbb{R}^n)}
	= O (t^{-\gamma_q-\frac{n}2+\frac12})
\end{equation}
as $t \to +\infty$.
The fourth term of \eqref{bsn1} is further separated to
\[
\begin{split}
	&\int_0^t \int_{\mathbb{R}^n} \biggl( \boldsymbol{a}\cdot\nabla G(t-s,x-y) - \sum_{2l+|\beta|=0}^{n-2} \frac{\partial_t^l \nabla^\beta (\boldsymbol{a}\cdot\nabla) G(t)}{l!\beta!} (-s)^l (-y)^\beta \biggr) u^2 (s,y) dyds\\
	&=
	J_{n-1} (t)+r_n^1 (t)
\end{split}
\]
for
\[
\begin{split}
	&J_{n-1} (t)
	= \int_0^t \int_{\mathbb{R}^n} \biggl( \boldsymbol{a}\cdot\nabla G(t-s,x-y) - \sum_{2l+|\beta|=0}^{n-2} \frac{\partial_t^l \nabla^\beta (\boldsymbol{a}\cdot\nabla) G(t)}{l!\beta!} (-s)^l (-y)^\beta \biggr)\\
	&\hspace{20mm} U_0^2 (s,y) dyds
\end{split}
\]
and
\[
\begin{split}
	 &r_{n}^1 (t)
	=
	\int_0^t \int_{\mathbb{R}^n} \biggl( \boldsymbol{a}\cdot\nabla G(t-s,x-y) - \sum_{2l+|\beta|=0}^{n-2} \frac{\partial_t^l \nabla^\beta (\boldsymbol{a}\cdot\nabla) G(t)}{l!\beta!} (-s)^l (-y)^\beta \biggr)\\
	&\hspace{20mm} (u^2-U_0^2) (s,y) dyds.
\end{split}
\]
The singularities of $U_0^2$ as $s \to +0$ are mitigated by Taylor theorem.
The last term can be further expanded to
\[
\begin{split}
	&r_n^1 (t)
	=
	\sum_{2l+|\beta|=n-1}  \frac{\partial_t^l \nabla^\beta (\boldsymbol{a}\cdot\nabla) G(t)}{l!\beta!} \int_0^t \int_{\mathbb{R}^n}
		(-s)^l (-y)^\beta (u^2 - U_0^2 - 2U_0 U_1) (s,y)
	dyds\\
	&+
	J_n (t)
	+ r_{n+1}^1 (t)
\end{split}
\]
for
\[
\begin{split}
	&J_n (t) = 2\int_0^t \int_{\mathbb{R}^n} \biggl( \boldsymbol{a}\cdot\nabla G(t-s,x-y) - \sum_{2l+|\beta|=0}^{n-1} \frac{\partial_t^l \nabla^\beta (\boldsymbol{a}\cdot\nabla) G(t)}{l!\beta!} (-s)^l (-y)^\beta \biggr)\\
	&\hspace{20mm} (U_0 U_1) (s,y) dyds
\end{split}
\]
and
\[
\begin{split}
	&r_{n+1}^1 (t)
	=
	\int_0^t \int_{\mathbb{R}^n} 
		\biggl( \boldsymbol{a}\cdot\nabla G(t-s,x-y) - \sum_{2l+|\beta|=0}^{n-1} \frac{\partial_t^l \nabla^\beta (\boldsymbol{a}\cdot\nabla) G(t)}{l!\beta!} (-s)^l (-y)^\beta \biggr)\\
		&\hspace{15mm} (u^2 - U_0^2 - 2U_0U_1) (s,y)
	dyds.
\end{split}
\]
Here $2U_0 U_1$ in the first part of $r_n^1$ is coming from \eqref{odd-van}.
Hence, there is no singularity here as $s \to +0$.
The singularity in $J_n$ is mitigated by Taylor theorem.
Similarly, $r_{n+1}^1$ has no singularity.
By repeating this procedure, we have that
\begin{equation}\label{bsn2}
\begin{split}
	&u(t) = \sum_{|\alpha| = 0}^{2n-3} \frac{\nabla^\alpha G(t)}{\alpha!} \int_{\mathbb{R}^n} (-y)^\alpha u_0 (y) dy\\
	&+ \sum_{2l+|\beta|=0}^{n-3} \frac{\partial_t^l \nabla^\beta (\boldsymbol{a}\cdot\nabla) G(t)}{l!\beta!} \int_0^t \int_{\mathbb{R}^n} (-s)^l (-y)^\beta
		u^2 (s,y)
	dyds\\
	&+ \sum_{2l+|\beta|=n-2}^{2n-4} \frac{\partial_t^l \nabla^\beta (\boldsymbol{a}\cdot\nabla) G(t)}{l!\beta!} \int_0^t \int_{\mathbb{R}^n} (-s)^l (-y)^\beta
		\biggl( u^2\\
		&\hspace{50mm} - \sum_{m_1+m_2=0}^{2l+|\beta|-n+2} U_{m_1} U_{m_2} \biggr) (s,y) 
	dyds\\
	&+ \sum_{m=n-1}^{2n-3} J_m (t) + r_{2n-2}^0 (t) + r_{2n-2}^1 (t)
\end{split}
\end{equation}
for
\[
\begin{split}
	&J_m (t) = \sum_{m_1 + m_2 = m-n+1} \int_0^t \int_{\mathbb{R}^n}
		\biggl( \boldsymbol{a}\cdot\nabla G(t-s,x-y)\\
		&\hspace{25mm} - \sum_{2l+|\beta|=0}^{m-1} \frac{\partial_t^l \nabla^\beta (\boldsymbol{a}\cdot\nabla) G(t,x)}{l!\beta!} (-s)^l (-y)^\beta \biggr)
		(U_{m_1}U_{m_2}) (s,y)
	dyds,
\end{split}
\]
and $r_{2n-2}^0$ given by \eqref{rm0}
and
\[
\begin{split}
	&r_{2n-2}^1 (t)
	=
	\int_0^t \int_{\mathbb{R}^n} 
		\biggl( \boldsymbol{a}\cdot\nabla G(t-s,x-y) - \sum_{2l+|\beta|=0}^{2n-4} \frac{\partial_t^l \nabla^\beta (\boldsymbol{a}\cdot\nabla) G(t)}{l!\beta!} (-s)^l (-y)^\beta \biggr)\\
		&\hspace{25mm} \biggl( u^2 - \sum_{m_1+m_2=0}^{n-2} U_{m_1} U_{m_2} \biggr) (s,y)
	dyds.
\end{split}
\]
Here the singularities of $U_{m_1} U_{m_2}$ on $J_m$ and $r_{2n-2}^1$ as $s \to +0$ are mitigated by Taylor theorem.
The error term $r_{2n-2}^1$ is treated in the similar way as in Section \ref{sect3dim}.
Precisely, Taylor theorem together with \eqref{decay-wt} and \eqref{EZexpK2} says that
$\| r_{2n-2}^1 (t) \|_{L^q (\mathbb{R}^n)} = O (t^{-\gamma_q-n+1} \log t)$ as $t \to +\infty$.
The coefficients on the second and third parts of \eqref{bsn2} are practically some polynomials $\mathcal{P}_{l\beta}$ of $t^{-1/2}$, respectively.
In fact, we derive $\mathcal{P}_{l\beta}$ for the second part as
\[
	\frac1{l!\beta!} \int_0^t \int_{\mathbb{R}^n} (-s)^l (-y)^\beta u^2 (s,y) dyds = \mathcal{P}_{l\beta} (t) + \eta_{l\beta} (t)
\]
for
\[
\begin{split}
	&\mathcal{P}_{l\beta} (t)
	= \frac1{l!\beta!} \int_0^\infty \int_{\mathbb{R}^n} (-s)^l (-y)^\beta u^2 (s,y) dyds\\
	&- \frac1{l!\beta!} \sum_{m_1 + m_2 = 0}^{n-2} \int_t^\infty s^{-\frac{n}2 - \frac{m_1}2 - \frac{m_2}2 + l + \frac{|\beta|}2} ds \int_{\mathbb{R}^n} (-1)^l (-y)^\beta (U_{m_1} U_{m_2}) (1,y)dy
\end{split}
\]
and
\[
	\eta_{l\beta} (t)
	= - \frac1{l!\beta!} \int_t^\infty \int_{\mathbb{R}^n} (-s)^l (-y)^\beta \biggl( u^2 - \sum_{m_1+m_2 = 0}^{n-2} U_{m_1} U_{m_2} \biggr) (s,y) dyds
\]
for $0 \le 2l+|\beta| \le n-3$.
For the third part of \eqref{bsn2}, we see that
\[
	\frac1{l!\beta!} \int_0^t \int_{\mathbb{R}^n} (-s)^l (-y)^\beta
		\biggl( u^2 - \sum_{m_1+m_2=0}^{2l+|\beta|-n+2} U_{m_1} U_{m_2} \biggr) (s,y) 
	dyds
	=
	 \mathcal{P}_{l\beta} (t) + \eta_{l\beta} (t)
\]
for
\[
\begin{split}
	&\mathcal{P}_{l\beta} (t)
	=
	\frac1{l!\beta!} \int_0^\infty \int_{\mathbb{R}^n} (-s)^l (-y)^\beta
		\biggl( u^2 - \sum_{m_1+m_2=0}^{2l+|\beta|-n+2} U_{m_1} U_{m_2} \biggr) (s,y) 
	dyds\\
	&- \frac1{l!\beta!} \sum_{m_1+m_2=2l+|\beta|-n+3}^{n-2} \int_t^\infty s^{-\frac{n}2-\frac{m_1}2-\frac{m_2}2+l+\frac{|\beta|}2} ds \int_{\mathbb{R}^n} (-1)^l (-y)^\beta\\
	&\hspace{60mm} (U_{m_1} U_{m_2}) (1,y) dy
\end{split}
\]
for $n-2 \le 2l+|\beta| \le 2n-5$, and
\[
	\frac1{l!\beta!} \int_0^t \int_{\mathbb{R}^n} (-s)^l (-y)^\beta
		\biggl( u^2 - \sum_{m_1+m_2=0}^{n-2} U_{m_1} U_{m_2} \biggr) (s,y) 
	dyds
	=
	 \mathcal{P}_{l\beta} (t) + \eta_{l\beta} (t)
\]
for
\[
\begin{split}
	&\mathcal{P}_{l\beta} (t)
	=
	\frac1{l!\beta!} \int_0^\infty \int_{\mathbb{R}^n} (-s)^l (-y)^\beta
		\biggl( u^2 - \sum_{m_1+m_2=0}^{n-2} U_{m_1} U_{m_2} \biggr) (s,y) 
	dyds
\end{split}
\]
for $2l+|\beta|=2n-4$.
The remaind terms $\eta_{l\beta}$ are same as above.
Here the constant terms of $\mathcal{P}_{l\beta}$ are sure integrable.
Indeed, since $U_{m_1} U_{m_2}$ for $m_1 + m_2 = 2l+|\beta|-n+2$ is coming from \eqref{odd-van}, we see that
\[
\begin{split}
	&\biggl| \int_0^\infty \int_{\mathbb{R}^n} (-s)^l (-y)^\beta
		\biggl( u^2 - \sum_{m_1+m_2=0}^{2l+|\beta|-n+2} U_{m_1} U_{m_2} \biggr) (s,y) 
	dyds \biggr|\\
	&\le C \int_0^\infty s^{-1/2} (1+s)^{-1} ds < \infty
\end{split}
\]
from \eqref{decay-wt} and \eqref{EZexpK2}.
%
Then, $u$ is finally written as
\[
\begin{split}
	&u(t) = \sum_{|\alpha| = 0}^{2n-3} \frac{\nabla^\alpha G(t)}{\alpha!} \int_{\mathbb{R}^n} (-y)^\alpha u_0 (y) dy
	+ \sum_{2l+|\beta|=0}^{2n-4} \mathcal{P}_{l\beta} (t) \partial_t^l \nabla^\beta (\boldsymbol{a}\cdot\nabla) G(t)
	+ \sum_{m=n-1}^{2n-3} J_m (t)\\
	&+ r_{2n-2}^0 (t) + r_{2n-2}^1 (t) + r_{2n-2}^2 (t)
\end{split}
\]
for
\[
	r_{2n-2}^2 (t)
	=  \sum_{2l+|\beta|=0}^{2n-4} \eta_{l\beta} (t) \partial_t^l \nabla^\beta (\boldsymbol{a}\cdot\nabla) G(t).
\]
From \eqref{decay-wt} and \eqref{EZexpK2}, we see
$
	|\eta_{l\beta} (t)|
	\le
	C t^{-n+\frac32+l+\frac{|\beta|}2}
$
and then $\| r_{2n-2}^2 (t) \|_{L^q (\mathbb{R}^3)} = O (t^{-\gamma_q-n+1})$ as $t \to +\infty$.
The terms of second part $\mathcal{P}_{l\beta} (t) \partial_t^l \nabla^\beta (\boldsymbol{a}\cdot\nabla) G(t)$ have clear scales.
The parabolic scales of $U_m$ provide that $\lambda^{n+m} J_m (\lambda^2 t, \lambda x) = J_m (t,x)$ for $\lambda > 0$ and then $J_m$ is a part of $U_m$.
Consequently, we decide concrete profiles $U_{n-1}, U_n,\ldots,$ $U_{2n-3}$ of higher-orders and never saw logarithms.

\appendix

\section{The normal convection-diffusion equation}
We introduce the nonlinear force term $\boldsymbol{a}\cdot\nabla (|u|u)$ in \eqref{cd} instead of $\boldsymbol{a}\cdot\nabla (u^2)$.
Since the parity is changed, we expect the asymptotic behavior of solutions to be completely different.
Even in this case, Escobedo and Zuazua \cite{EZ} derived the corresponding asymptotic profiles $U_m$ for $0 \le m \le n-2$ and proved \eqref{EZexp}.
Unfortunatelly, in this case, our renormalization does not work.
For example, the error term $r_{2n-2}^1$ in Section \ref{sect-oddim} changes to
\[
\begin{split}
	&r_{2n-2}^1 (t)
	=
	\int_0^t \int_{\mathbb{R}^n} 
		\biggl( \boldsymbol{a}\cdot\nabla G(t-s,x-y) - \sum_{2l+|\beta|=0}^{2n-4} \frac{\partial_t^l \nabla^\beta (\boldsymbol{a}\cdot\nabla) G(t)}{l!\beta!} (-s)^l (-y)^\beta \biggr)\\
		&\hspace{25mm} \biggl( |u|u - \sum_{m_1+m_2=0}^{n-2} |U_{m_1}| U_{m_2} \biggr) (s,y)
	dyds
\end{split}
\]
in this case.
To estimate this quantity by using \eqref{EZexp}, we should guarantee that
\[
	\biggl| |u| - \sum_{m = 0}^{n-2} |U_{m}| \biggr|
	\le
	C \biggl| u - \sum_{m = 0}^{n-2} U_{m} \biggr|.
\]
Needless to say, such an inequality is false.

\end{document}